\documentclass[12pt,a4paper]{article}
\usepackage[english]{babel}
\usepackage[utf8]{inputenc}
\usepackage{amsmath}
\usepackage{amsfonts}
\usepackage{amssymb}
\usepackage{graphicx}
\usepackage{listings}
\usepackage{listingsutf8}
\usepackage{tikz}
\usepackage{cite}
\usepackage{algorithm2e}
\usepackage{array}
\usepackage{geometry}
\usepackage{makecell}
\usepackage{colortbl}
\usepackage{float}
\DeclareMathOperator{\rank}{rk}
\DeclareMathOperator{\trace}{Tr}

\begin{document}

\allowdisplaybreaks

\setcounter{tocdepth}{2}

\nocite{*}
\newcounter{currentcounter}

\newcounter{corollarycounter}

\newcounter{example}
\newcommand\exc{\theexample}

\newenvironment{general}[1]{\refstepcounter{currentcounter} \bigskip \noindent \textbf{#1 \thecurrentcounter :}\itshape}{\normalfont \medskip}

\newenvironment{generalbis}[2]{\refstepcounter{currentcounter} \bigskip \noindent \textbf{#1 \thesubsection.\thecurrentcounter ~#2 :}\itshape}{\normalfont \medskip}

\newenvironment{proof}{\noindent \textbf{Proof :} \newline \noindent \hspace*{0.2cm}}{\hspace*{\fill}$\square$ \bigskip}

\newenvironment{proofbis}[1]{\noindent \textbf{#1:} \newline \noindent \hspace*{0.2cm}}{\hspace*{\fill}$\square$ \bigskip}

\newenvironment{famous}[1]{\medskip \noindent \textbf{#1:} \newline \noindent  \itshape}{\normalfont \medskip}

\newenvironment{mycor}{\refstepcounter{corollarycounter} \medskip \noindent \textbf{Corollary \thecorollarycounter :}  \itshape}{\normalfont \medskip}

\newenvironment{court}[1]{\refstepcounter{currentcounter} \smallskip \noindent \textbf{#1 \thesubsection.\thecurrentcounter :} \newline \noindent \itshape}{\normalfont \smallskip}

\newcolumntype{C}{>{$}c<{$}}
\newcommand\toline{\smallskip \newline}
\newcommand\refp[1]{(\ref{#1})}
\renewcommand\mod{\mathrm{~mod~}}
\newcommand\pregcd{\mathrm{gcd}}
\renewcommand\gcd[2]{\pregcd(#1,#2)}
\newcommand\prelcm{\mathrm{lcm}}
\newcommand\lcm[2]{\prelcm(#1,#2)}
\newcommand\naturaliso{\cong}
\newcommand\iso{\simeq}
\newcommand\cross[1]{#1^{\times}}
\newcommand\crosslong[1]{\cross{(#1)}}
\newcommand\poldegree[1]{\mathrm{deg}(#1)}
\newcommand\cardinal[1]{\# \left(#1\right)}
\newcommand\cardinalshort[1]{\##1}
\renewcommand\det{\mathrm{det}}
\newcommand\indicator[1]{[#1]}

\newcommand\sign{\mathrm{sign}}
\newcommand\signature{\mathrm{sgn}}

\newcommand\bb[1]{\mathbb{#1}}
\newcommand\kk{\bb{K}}
\newcommand\ok{\mathcal{O}_{\kk}}
\newcommand\zz{\bb{Z}}
\newcommand\qq{\bb{Q}}
\newcommand\rr{\bb{R}}
\newcommand\cc{\bb{C}}
\newcommand\ff{\bb{F}}
\newcommand\of{\mathcal{O}_{\ff}}
\newcommand\goth[1]{\mathfrak{#1}}
\newcommand\zsz[1]{\zz/#1\zz}
\newcommand\zszcross[1]{\zsz{#1}^{\times}}
\newcommand\dualsimple[1]{#1^{\vee}}
\newcommand\dual[1]{(#1)^{\vee}}

\newcommand\dx[1]{\mathrm{d}#1}
\newcommand\ddx[1]{\frac{\mathrm{d}}{\dx{#1}}}
\newcommand\dkdxk[2]{\frac{\mathrm{d}^{#2}}{\dx{#1}^{#2}}}
\newcommand\partialx[1]{\partial #1}
\newcommand\partialdx[1]{\frac{\partial}{\partialx{#1}}}
\newcommand\partialkdxk[2]{\frac{\partial^{#2}}{\partialx{#1}^{#2}}}
\newcommand\integral[3]{\int_{#1}{#2}\dx{#3}}
\newcommand\integralsimple[1]{\integrale{G}{#1}{\lambda}}
\newcommand\sprod[2]{\langle #1, #2 \rangle}

\newcommand\limi[1]{\lim_{#1 \to + \infty}}

\newcommand\mn[2]{\mathrm{M}_{#1}(#2)}
\newcommand\mnz[1]{\mn{#1}{\zz}}
\newcommand\sln[2]{\mathrm{SL}_{#1}(#2)}
\newcommand\slnz[1]{\sln{#1}{\zz}}
\newcommand\gln[2]{\mathrm{GL}_{#1}(#2)}
\newcommand\glnz[1]{\gln{#1}{\zz}}
\newcommand\resc{\mathrm{resc}}
\newcommand\com{\mathrm{com}}

\newcommand\bars[1]{\underline{#1}}
\newcommand\taubar{\bars{\tau}}
\newcommand\qbar{\bars{q}}
\newcommand\sigmabar{\bars{\sigma}}
\newcommand\rbar{\bars{r}}
\newcommand\nbar{\bars{n}}
\newcommand\mbar{\bars{m}}
\newcommand\mubar{\bars{\mu}}
\newcommand\xbar{\bars{x}}
\newcommand\Xbar{\bars{X}}
\newcommand\abar{\bars{a}}
\newcommand\bbar{\bars{b}}
\newcommand\alphabar{\bars{\alpha}}
\newcommand\betabar{\bars{\beta}}
\newcommand\omegabar{\bars{\omega}}
\newcommand\taubarsj[1]{\taubar^{-}(#1)}
\newcommand\onebar{\bars{1}}

\newcommand\hh{\bb{H}}
\newenvironment{psmallmatrix}
  {\left(\begin{smallmatrix}}
  {\end{smallmatrix}\right)}
  
\newcommand\declarefunction[5]{#1 := \begin{cases}
\hfill #2 \hfill & \to ~ #3 \\
\hfill #4 \hfill & \to ~ #5 \\
\end{cases}}

\newcommand\sume{\sideset{}{_e}\sum}

\newcommand\tendstowhen[1]{\xrightarrow[#1]{}}

\newcommand\pending{$\bigskip \newline \blacktriangle \blacktriangle \blacktriangle \blacktriangle \blacktriangle \blacktriangle \blacktriangle \blacktriangle \bigskip \newline \hfill$}

\newcommand\pendingref{\textbf{[?]}}

\newcommand\opc[1]{\mathcal{O}^{+,\times}_{#1}}
\newcommand\opck{\opc{\kk}}
\newcommand\opcf{\opc{\goth{f}}}
\newcommand\eps{\varepsilon}
\newcommand\units{\ok^{\times}}
\newcommand\unitsf{\of^{\times}}
\newcommand\norm[1]{\mathcal{N}(#1)}

\newcommand\ta{\tilde{a}}
\newcommand\tabar{\tilde{\abar}}
\newcommand\talpha{\tilde{\alpha}}
\newcommand\talphabar{\tilde{\alphabar}}
\newcommand\ts{\tilde{s}}
\newcommand\ttt{\tilde{t}}
\newcommand\tlambda{\tilde{\lambda}}
\newcommand\tc{\tilde{c}}
\newcommand\tcone{\tc_1}
\newcommand\tctwo{\tc_2}
\newcommand\te{\tilde{e}}
\newcommand\teone{\te_1}
\newcommand\tetwo{\te_2}
\newcommand\tmu{\tilde{\mu}}
\newcommand\tA{\tilde{A}}

\newcommand\spvgamma{\Gamma_{\goth{f}, \goth{b}, \goth{a}}(\eps; h)}
\newcommand\spvgammasign{\Gamma_{\goth{f}, \goth{b}, \goth{a}}^{\pm}(\eps; h)}
\newcommand\spvgr{G_{r, \goth{f}, \goth{b}, \goth{a}}^{\pm}(u_1, \dots, u_r; h)}
\newcommand\spvgrrho{G_{r, \goth{f}, \goth{b}, \goth{a}}^{\pm}([\eps_{\rho(1)}|\dots|\eps_{\rho(r)}]; h_\rho)}
\newcommand\spvgrcomplete{I_{r,\goth{f}, \goth{b}, \goth{a}}(\eps_1, \dots, \eps_r ; \bars{h}, \bars{\mu}, \bars{\nu})}

\newcommand\spvgammac{\Gamma_{\goth{f}, \goth{b}, \goth{a}}(\eps; h, \sigma_{\cc})}
\newcommand\spvgrc{G_{r, \goth{f}, \goth{b}, \goth{a}}^{\pm}([\eps_{\rho(1)}|\dots|\eps_{\rho(r)}]; h_\rho, \sigma_{\cc})}
\newcommand\spvgrrhoc{G_{r, \goth{f}, \goth{b}, \goth{a}}^{\pm}([\eps_{\rho(1)}|\dots|\eps_{\rho(r)}]; h_\rho)}
\newcommand\spvgrcompletec{I_{r,\goth{f}, \goth{b}, \goth{a}}(\eps_1, \dots, \eps_r ; \bars{h}, \bars{\mu}, \bars{\nu}, \sigma_{\cc})}

\newcommand\ula{u_{L, \goth{a}}}
\newcommand\ulah{u_{L, \goth{a}, \bars{h}}}
    
\newcommand\hombase{\mathrm{Hom}}
\newcommand\myhom[3]{\hombase_{#1}(#2, #3)}    
\newcommand\homlong[3]{\hombase_{#1}\left(#2, #3\right)}   
\newcommand\homlz{\myhom{\zz}{L}{\zz}}
\newcommand\homlprimez{\myhom{\zz}{L'}{\zz}}
\newcommand\homlc{\myhom{\zz}{L}{\cc}}
\newcommand\homlambdaz{\myhom{\zz}{\Lambda}{\zz}}
\newcommand\homlambdac{\myhom{\zz}{\Lambda}{\cc}}
\newcommand\homvq{\myhom{\qq}{V}{\qq}}
\newcommand\zexc{z}
  
\newcommand\zfone{\mathcal{Z}_{\goth{f}}^1}
  
\newcommand\fracpart[1]{\left\{#1\right\}}
\newcommand\entirepart[1]{\left\lfloor#1\right\rfloor}  
\newcommand\shortsetseparator[1][]{#1|}
\newcommand\setseparator[1][]{~\shortsetseparator[#1]~}
\newcommand\coeff{\mathrm{coeff}}
\newcommand\normalbone[1]{\left(\left(#1\right)\right)}
\newcommand\tensor{\otimes}
\newcommand{\at}[2][]{#1|_{#2}}
\newcommand\omitvar[1]{\widehat{#1}}
\newcommand\kronecker{\delta}
\newcommand\plgt{x}
\newcommand\cohomd{\partial}
\newcommand\cohomdx{\partial^{\times}}
\newcommand\compset[1]{#1^c}
\newcommand\signdet{\sign\,\det}

\newcommand\copen{c^{\circ}}
\newcommand\cclosed{c}
\newcommand\cdual{c^{\vee}}
\newcommand\copendual{c^{\vee, \circ}}
\newcommand\dirac{\delta}

\newcommand\myspan{\mathrm{Span}}
\newcommand\prespanconvex[1]{\mathcal{C}(#1)}
\newcommand\spanconvex{\prespanconvex{V}}
\newcommand\prespancones[1]{\mathcal{K}(#1)}
\newcommand\spancones{\prespancones{V}}
\newcommand\spanconesrr{\prespancones{V_{\rr}}}
\newcommand\spanconesk[1]{\mathcal{K}^{#1}(V)}
\newcommand\spanqcones{\mathcal{K}_{\qq}(V)}
\newcommand\spanqconesrr{\mathcal{K}_{\qq}(V_{\rr})}
\newcommand\prespanwedges[1]{\mathcal{L}(#1)}
\newcommand\spanwedges{\prespanwedges{V}}
\newcommand\spanqwedges{\mathcal{L}_{\qq}(V)}
\newcommand\spanqwedgesrr{\mathcal{L}_{\qq}(V_{\rr})}

\newcommand\coefficient[3]{\mathrm{coeff}[#2^{#3}]\left(#1\right)}
\newcommand\badposition{\mathrm{(BP)}}
\newcommand\symbolparagraph{\S}

\newcommand\alphajk[2]{\alpha^{(#1)}_{#2}}
\newcommand\djk[2]{d^{(#1)}_{#2}}
\newcommand\vjkl{v^{(j)}_{k,k'}}
\newcommand\yj[1]{y^{(#1)}}
\newcommand\gammaj[1]{\gamma^{(#1)}}
\newcommand\yjk[2]{y^{(#1)}_{#2}}
\newcommand\ujk[2]{u^{(#1)}_{#2}}
\newcommand\Ujk[2]{U^{(#1)}_{#2}}

\newcommand\disjointunion{\sqcup}
\newcommand\dkxs{d_k(x, s)}

\newcommand\hexp{h_e}

\newcommand\smoothedgr{G_{n-2, a_1, \dots, a_{n-1}}(v)(w, x, L, L')}
\newcommand\smoothedgromitj{G_{n-2, a_1, \dots, \omitvar{a_j}, \dots, a_{n}}(v)(w, x, L, L')}
\newcommand\smoothedbn{B_{n, a_1, \dots, a_n}(v)(w,x, L, L')}
\newcommand\congruencegrouphnn{\Gamma_0(N, n)}
\newcommand\generaldenomclasscohom[1]{\mathcal{D}(#1)}
\newcommand\denomclasscohom{\generaldenomclasscohom{N, n}}

\newcommand\cun[1]{c^1_{#1}}
\newcommand\cdeux[1]{c^2_{#1}}
\newcommand\Cun[1]{C^{1}_{#1}}
\newcommand\Cdeux[1]{C^{2}_{#1}}
\newcommand\fun{f^1}
\newcommand\fdeux{f^2}
\newcommand\mathcalCpm[2]{\mathcal{C}^{#1}_{#2}}
\newcommand\mathcalCplus[1]{\mathcal{C}^{+}_{#1}}
\newcommand\mathcalCmoins[1]{\mathcal{C}^{-}_{#1}}
\newcommand\good{\textit{good }}
\newcommand\hj[2]{H^{#1}_{#2}}
\newcommand\hjplus[1]{\hj{+}{#1}}
\newcommand\hjminus[1]{\hj{-}{#1}}

\newcommand\lemmasignstranspose{lemma \ref{lemmasignrelations}, (ii) }
\newcommand\lemmasignstoprows{lemma \ref{lemmasignrelations}, (i) }
\newcommand\lemmasignsujk{lemma \ref{lemmasignrelations}, (iii) }

\newcommand\wellplaced{well placed }

\begin{center}
\large \noindent Geometric families of multiple elliptic Gamma functions and arithmetic applications, II

\medskip \noindent \textit{Action of congruence subgroups in} $\slnz{n}$
\end{center}
\medskip
\begin{center}
Pierre L. L. Morain\footnotemark[1]\footnotetext[1]{Sorbonne Université and Université Paris Cité, CNRS, INRIA, IMJ-PRG, F-75005 Paris, France. This PhD work is funded by the École polytechnique,  Palaiseau, France.}
\end{center}
\medskip
\begin{center}
\noindent \textbf{Abstract:}
\end{center}
This is the second paper in a series where we study arithmetic applications of the multiple elliptic Gamma functions originated in mathematical physics. In the first article in this series we defined geometric families of these functions and proved that these families satisfied coboundary relations involving an attached collection of Bernoulli rational functions. The main purpose of the present paper is to show that smoothed versions of our geometric elliptic Gamma functions give rise to partial modular symbols for congruence subgroups of $\mathrm{SL}_{n}(\mathbb{Z})$ for $n \geq 2$ which restrict to $(n-2)$-cocycles on tori in $\mathrm{SL}_{n}(\mathbb{Z})$ coming from groups of totally positive units in number fields. To achieve this, we show that the associated smoothed Bernoulli rational functions reduce to smoothed higher Dedekind sums with uniformly bounded denominators.

\normalsize \bigskip \bigskip

\noindent \textbf{Acknowledgments:} 
This article was written as part of an on-going PhD work and the author would like to thank his advisors Pierre Charollois and Antonin Guilloux at Sorbonne Université for their guidance and their helpful comments, as well as their support throughout this journey. The author is also grateful to Mateo Crabit Nicolau, Mart\'i Roset Julià and Peter Xu for helpful discussions.
\newpage

\tableofcontents

\section{Introduction}

In this series of papers we are interested in multivariate meromorphic functions built on the model of a $\theta$ function used in the theory of \textit{Complex Multiplication} which is defined on $\cc \times \hh$ by:
$$\theta(z, \tau) = \prod_{n \geq 0} \left(1- e^{-2i\pi z} e^{2i \pi (n+1) \tau}\right) \left(1-e^{2i\pi z}e^{2i\pi n \tau}\right)$$
where $\hh$ denotes the upper half-plane. This function enjoys modular properties under an action of $\slnz{2}$. Indeed, there is a rational function $P_{2} : \slnz{2} \to \qq(z, \tau)$ such that for all $\gamma = \begin{psmallmatrix} a & b \\ c & d \end{psmallmatrix}$:
\begin{equation}\label{dedekindrademacher}
\theta\left(\frac{z}{c\tau +d}, \frac{a\tau + b}{c\tau + d}\right) = \theta(z, \tau). e^{2i\pi P_{2, \gamma}(z, \tau)}.
\end{equation}
The rational function $P_{2, \gamma}$ may be expressed in terms of the classic Dedekind sums defined for $c > 0$ by:
\begin{equation}\label{dedekindsums}
s(c, d) = \sum_{k = 1}^{c-1} b_1\left(\frac{k}{c}\right) b_1\left(\frac{kd}{c}\right)
\end{equation}
where $b_1 : t \to t - \lfloor t \rfloor - 1/2$ is the classic periodic version of the Bernoulli polynomial $B_1(t) = t - 1/2$. The explicit formula for the rational function $P_{2, \gamma}$ is given by:
$$P_{2, \gamma}(z, \tau) = \begin{cases} 0 & \text{ if } c = 0 \text{ and } a = 1 \\ z- \frac12 & \text{ if } c = 0 \text{ and } a = -1 \\ \frac{z^2c^2 +zc + 1/6}{2c(c\tau + d)} - \frac{z}{2}+\frac{c\tau + d}{12c} - \sign(c)(s(|c|, d)+\frac{1}{4}) & \text{ if } c \neq 0\end{cases} $$
where $\sign(c) = c/|c|$ when $c \neq 0$ (this can be derived from [\hspace{1sp}\cite{Siegel}, Proposition 4] in conjunction with [\hspace{1sp}\cite{Rademacher}, Formula (3.26)]). It is classical to perform a smoothing operation on the $\theta$ function to obtain modular units of one of the two following shapes:
$$\frac{\theta(z, \tau)^{N^2}}{\theta(Nz, \tau)}~~\textit{ or } ~~ \frac{\theta(z, \tau)^N}{\theta(Nz, N\tau)}.$$
In the first case, by adding a small exponential prefactor we obtain the basic ingredient for Siegel units:
$$ _N\theta(z, \tau) = e^{2i\pi\left(\frac{N^2-1}{12}\tau + \frac{N - N^2}{2}\left(z-\frac{1}{2}\right)\right)}\frac{\theta(z, \tau)^{N^2}}{\theta(Nz, \tau)} $$
when $z \in \qq$, $\tau$ is an imaginary quadratic number and $N$ is a rational prime which is inert in the field $\qq(\tau)$ (see [\hspace{1sp}\cite{Kato}, Proposition 1.3]). The second expression gives smoothed versions of Robert's elliptic units \cite{Robert} in the setting where $N$ is a rational prime which splits in $\qq(\tau)$. In this series we focus on the second type of smoothing and if we set:
$$\theta^{(N)}(z, \tau) = \frac{\theta(z, \tau)^N}{\theta(Nz, N\tau)}$$
then \refp{dedekindrademacher} gives for any $N \mid c$:
\begin{equation}\label{dedekindrademachersmooth}\theta^{(N)}\left(\frac{z}{c\tau +d}, \frac{a\tau + b}{c\tau + d}\right) = \theta^{(N)}(z, \tau). e^{2i\pi P_{2, \gamma}^{(N)}(z, \tau)}
\end{equation}
where 
$$P_{2, \gamma}^{(N)}(z, \tau) = \begin{cases} 0 & \text{ if } c = 0 \text{ and } d = 1 \\ \frac{1-N}{2} & \text{ if } c = 0 \text{ and } d = -1  \\ \sign(c)\left(s\left(\frac{|c|}{N}, d\right) - Ns(|c|, d) + \frac{1-N}{4}\right) & \text{ if } c \neq 0.\end{cases} $$
Crucially, the smoothed modularity defect $P_{2, \gamma}^{(N)}$ depends only on the matrix $\gamma \in \Gamma_0(N) = \{ \begin{psmallmatrix} a & b \\ c & d\end{psmallmatrix} \setseparator ~~N \mid c \}$ and not on $z$ and $\tau$. Moreover, it satisfies an integrality property as the values of $P_{2, \gamma}^{(N)}$ lie in $\frac{1}{12}\zz$. The application $(\gamma \to P_{2, \gamma}^{(N)})$ is essentially a smoothed version of the classic Dedekind-Rademacher function. This function is defined on $\slnz{2}$ by:
$$\phi_{DR}\begin{pmatrix} a & b \\ c & d \end{pmatrix} := \begin{cases} \frac{b}{d} & \text{ if } c = 0 \\ \frac{a + b}{c} - 12.\sign(c) s(|c|,d) & \text{ if } c \neq 0\end{cases} $$
and it takes values in $\zz$ (see \cite{Rademacher}). The function $\phi_{DR}$ is almost a group morphism as for any $\gamma, \gamma' \in \slnz{2}$:
$$\phi_{DR}(\gamma'') = \phi_{DR}(\gamma) + \phi_{DR}(\gamma') -3.\sign(cc'c'') $$
where $\gamma'' = \begin{psmallmatrix} a'' & b'' \\ c'' & d'' \end{psmallmatrix} = \begin{psmallmatrix} a' & b' \\ c' & d' \end{psmallmatrix} \begin{psmallmatrix} a & b \\ c & d \end{psmallmatrix}= \gamma'\gamma$. 
The level $N$ function $\Psi_N$ defined on $\Gamma_0(N)$ by 
$$\Psi_N \begin{pmatrix} a & b \\ c & d \end{pmatrix} =\phi_{DR}\begin{pmatrix} a & bN \\ \frac{c}{N} & d \end{pmatrix} - \phi_{DR}\begin{pmatrix} a & b \\ c & d \end{pmatrix} $$
(see \cite{Rademacher}, \cite{Mazur}) then defines a $1$-cocycle in $H^1(\Gamma_0(N), \zz)$. The Dedekind-Rademacher function and its level $N$ avatars have played a key role in the construction of conjectural elliptic units above real quadratic fields by Darmon and Dasgupta \cite{DarmonDasgupta} and they are also connected to Gross-Stark units above real quadratic fields (see \cite{DarmonPozziVonk}). In our context the application $(\gamma \to P_{2, \gamma}^{(N)})$ is given by a different smoothed version of the Dedekind-Rademacher function:
$$P_{2, \gamma}^{(N)} = \frac{1}{12}\left(N\phi_{DR}\begin{pmatrix} a & b \\ c & d \end{pmatrix} -\phi_{DR}\begin{pmatrix} a & bN \\ \frac{c}{N} & d \end{pmatrix}\right) + \begin{cases} 0 & \text{ if } c = 0 \text{ and } d = 1 \\ \frac{1-N}{2} & \text{ if } c = 0 \text{ and } d = -1 \\ \sign(c)\frac{1-N}{4} & \text{ if } c \neq 0 \end{cases}  $$
and its reduction mod $\zz$ defines a $1$-cocycle on $\Gamma_0(N)$ with values in $\frac{1}{12}\zz/\zz$. 

Our goal throughout this series is to generalise the construction of elliptic units to general number fields with exactly one complex place using higher degree analogues of the $\theta$ function satisfying similar transformation properties for higher degree special linear groups $\slnz{n}$, $n \geq 2$. To achieve this, we generalise the construction of conjectural elliptic units above complex cubic fields (that is $n =3$) by Bergeron, Charollois and Garc\'ia to higher degree number fields with exactly one complex place. The conjectural elliptic units in \cite{BCG} are given by specific evaluations of the elliptic Gamma function of Ruijsenaars \cite{Ruijsenaars}. This function, which was studied at length by Felder and Varchenko (see \cite{FV}), is a meromorphic function on $\cc \times \hh^2$ defined by:
$$\Gamma(z, \tau, \sigma) = \prod_{m,n \,\geq 0}\left(\frac{1-\exp(2i\pi((m+1)\tau +(n+1)\sigma-z))}{1-\exp(2i\pi(m\tau+n\sigma+z))} \right).$$
It satisfies modular transformation properties which allow for the construction of a $1$-cocycle for $\slnz{3}$. In the early 2000s Nishizawa \cite{Nishizawa} introduced a whole hierarchy of multiple elliptic Gamma functions to generalise both the theta and the elliptic Gamma function. They are meromorphic functions on $\cc \times \hh^{r+1}$ defined by:
$$G_r(z, \tau_0, \dots, \tau_r) = \prod_{m_0, \dots, m_r \,\geq 0} \left(1 - e^{2i\pi(-z + \sum_{j = 0}^r (m_j+1)\tau_j)}\right)  \left(1- e^{2i\pi (z + \sum_{j = 0}^r m_j\tau_j)}\right)^{(-1)^r}$$
for all $r \in \zz_{\geq 0}$ with the identification $\theta = G_0$ and $\Gamma = G_1$. They also satisfy similar modular transformation properties for $\slnz{r+2}$ as proven by Narukawa \cite{Narukawa}, involving some generalised Bernoulli polynomials as their modularity defect. The higher elliptic units we construct in \cite{thirdpaper} are given by specific evaluations of $G_{n-2}$ functions at points of a degree $n$ number field $\kk$ with exactly one complex place, thus giving a conjectural analytic description of the abelian extensions of $\kk$ in the spirit of Hilbert's 12th problem. 

In the first article in this series \cite{firstpaper} we introduced two geometric families of multiple elliptic Gamma functions and their attached Bernoulli rational functions. This construction generalises the construction of geometric families of the elliptic Gamma function carried out by Felder, Henriques, Rossi and Zhu \cite{FDuke} for the entire hierarchy of multiple elliptic Gamma functions. These geometric families consist of functions defined on lattices of rank $n \geq 2$ which enjoy modular properties for special linear groups $\slnz{n}$. More precisely, if $L$ is a lattice of rank $n$ in a $\qq$-vector space $V$ then the multiple elliptic Gamma function attached to $n-1$ linearly independent linear forms $a_1, \dots, a_{n-1} \in \homlz$ is a function $G_{n-2, a_1, \dots, a_{n-1}} : V/L \times \cc \times \homlc \to \cc$ defined by
\begin{multline*}
G_{n-2, a_1, \dots, a_{n-1}}(v)(w,x)^{(-1)^D} = \prod_{\delta \in (v+L) \cap C^{-}(\abar, x)/\zz\gamma} \left(1-e^{-2i\pi\left(\frac{w+\plgt(\delta)}{\plgt(\gamma)}\right)}\right) \\ \times \prod_{\delta \in (v+L) \cap C^{+}(\abar, x)/\zz\gamma} \left(1-e^{2i\pi\left(\frac{w+\plgt(\delta)}{\plgt(\gamma)}\right)}\right)^{(-1)^n}
\end{multline*}
for some explicit cones $C^{\pm}(\abar, x)$ in $V$, an explicit vector $\gamma \in L$ and an explicit integer $D \in \zz$ (see section \ref{sectionrecall} for a precise definition). As a function of $w, x \in \cc \times \homlc \simeq \cc \times \cc^n$ it is meromorphic. The associated Bernoulli rational function $B_{n, a_1, \dots, a_n}(v) \in \qq[w](x)$ is defined as the coefficient of $t^0$ in the expansion of the function
$$ \signdet(a_1,\dots, a_n)\times \sum_{\delta \in (v + L) \cap C} e^{wt + x(\delta)t}$$
where the cone $C$ is defined by $C = \{ \delta \in V \setseparator \forall\, 1 \leq j \leq n, a_j(\delta) \geq 0\}$ (see \refp{defbernoulliexplicit} for a precise definition). These families of functions are equivariant under the action of $\slnz{n}$ on the lattice $L$ as for all $g \in \slnz{n}$ they satisfy:
\begin{align}
G_{n-2,g \cdot a_1, \dots, g \cdot a_{n-1}}(g \cdot v)(w, g\cdot x) &= G_{n-2, a_1, \dots, a_{n-1}}(v)(w, x)\label{equivgr} \\
B_{n, g\cdot a_1, \dots, g\cdot a_n}(g \cdot v)(w, g\cdot x) &= B_{n, a_1, \dots, a_n}(v)(w, x).\label{equivbn}
\end{align}
More deeply, they satisfy partial coboundary and cocycle relations. First, if $a_1, \dots, a_n$ are linearly independent then:
\begin{equation}\label{modularproperty}
\prod_{j = 1}^{n} G_{n-2, a_1, \dots, \omitvar{a_j}, \dots, a_n}(v)(w, \plgt)^{(-1)^{j+1}} = \exp(2i\pi B_{n, a_1, \dots, a_n}(v)(w,\plgt))
\end{equation}
(see [\hspace{1sp}\cite{firstpaper}, Theorem 1]) which is a higher degree analogue of \refp{dedekindrademacher}. Secondly, for most configurations of $n+1$ linear forms $a_0, \dots, a_{n}$ the following cocycle relation holds:
\begin{equation}\label{cocyclebernoulli}
\sum_{j = 0}^n (-1)^j B_{n, a_0, \dots, \omitvar{a_j}, \dots, a_n}(v)(w,x) = 0
\end{equation}
(see [\hspace{1sp}\cite{firstpaper}, Corollary 1]). The main goal of this article is to prove that the modular property \refp{modularproperty} becomes a cocycle relation for the family of $G_{n-2, a_1, \dots, a_{n-1}}$ functions when a certain classic smoothing operation is performed, thus generalising \refp{dedekindrademachersmooth}. To state our main theorem we now introduce some notations relative to this smoothing operation. For the rest of this paper, we fix a rank $n \geq 2$ lattice $L$ with a $\zz$-basis $B = [e_1, \dots, e_n]$ as well as an integer $N \geq 2$. The lattice $L'$ generated by the $\zz$-basis $B' = [Ne_1, e_2, \dots, e_n]$ is called the smoothing lattice. We denote by $\Lambda$ (resp. $\Lambda'$) the dual space $\homlz$ (resp. $\homlprimez$) of $L$ (resp. $L'$) and define 
$$\Lambda_N = \{a \in \Lambda \setseparator a_{|L'} \text{ is primitive in } \Lambda' \} $$
where we recall that an element $a \in \Lambda$ is primitive if $a/d \in \Lambda$ for some integer $d$ implies $d = \pm 1$. The set $\Lambda_N$ is naturally endowed with an action of the following congruence subgroup in $\slnz{n}$:
\begin{equation}\label{definitionmodulargroup}
\congruencegrouphnn = \left\{ g \in \slnz{n} \setseparator[\Big] g \equiv \begin{pmatrix} \ast & \ast &\dots & \ast \\ 0 & \ast & \dots & \ast \\ \vdots & \vdots & \vdots & \vdots \\ 0 & \ast & \dots & \ast \end{pmatrix} \mod N \right\},
\end{equation}
the action being given by right multiplication by the inverse as $g \cdot a = a \times g^{-1}$. When the linear forms $a_1, \dots, a_n$ belong to $\Lambda_N$ it makes sense to introduce the smoothed functions:
\begin{equation}\label{defsmoothedgr}\smoothedgr = \frac{G_{n-2, a_1, \dots, a_{n-1}}(v)(w, x, L')^N}{G_{n-2, a_1, \dots, a_{n-1}}(v)(w, x, L)}, \end{equation}
\begin{equation}\label{defsmoothedbn}\smoothedbn = NB_{n, a_1, \dots, a_n}(v)(w,x, L') - B_{n, a_1, \dots, a_n}(v)(w,x,L).\end{equation}
We deduce immediately from \refp{equivgr} and \refp{equivbn} that these functions are equivariant under the action of $\congruencegrouphnn$ as for all $g \in \congruencegrouphnn$:
\begin{align}
G_{n-2,g \cdot a_1, \dots, g \cdot a_{n-1}}(g \cdot v)(w, g\cdot x, L, L') &= G_{n-2, a_1, \dots, a_{n-1}}(v)(w, x, L, L')\label{equivgrsmooth} \\
B_{n, g\cdot a_1, \dots, g\cdot a_n}(g \cdot v)(w, g\cdot x, L, L') &= B_{n, a_1, \dots, a_n}(v)(w, x, L, L').\label{equivbnsmooth}
\end{align}
It also follows at once from formula \refp{modularproperty} that these smoothed functions satisfy the coboundary relation:
\begin{equation}\label{smoothedmodularproperty}
\prod_{j = 1}^{n} \smoothedgromitj^{(-1)^{j+1}} = \exp(2i\pi \smoothedbn)
\end{equation}
for linearly independent $a_1, \dots, a_n \in \Lambda_N$. 

The goal of this paper is to derive a cocycle relation for the smoothed functions $\smoothedgr$ from relation \refp{smoothedmodularproperty}, turning the function 
$$a_1, \dots, a_{n-1} \to \smoothedgr$$ 
into a partial modular symbol for $\congruencegrouphnn$. Our main theorem is expressed under two conditions on the positions of the linear forms $a_1, \dots, a_n$ in $\Lambda$ relative to the smoothing lattice $L'$. First, we shall say that $a_1, \dots, a_n$ are \wellplaced (in $\dualsimple{V}$) if either $\rank(a_1, \dots, a_n) \neq n-1$ or if $\rank(a_1, \dots, a_n) = n-1$ and $0$ is not a barycenter of $a_1, \dots, a_n$ in $\dualsimple{V}$ (see Definition \ref{definitiongoodposition} for more details on this hypothesis). The second condition concerns the position of the lattice $L'$ in the case where $\rank(a_1, \dots, a_n) = n$. When $a_1, \dots, a_n$ are linearly independent, there are unique primitive elements $\alpha_1, \dots, \alpha_n$ in $L$ such that 
$$a_j(\alpha_j) > 0 ~~ \text{ and } ~~ a_j(\alpha_k) = 0, \forall\,k \neq j.$$
We shall say in the spirit of \cite{DasguptaShintani} that the index $N$ smoothing lattice $L'$ is \good for the linear forms $a_1, \dots, a_n$ if and only if $a_1, \dots, a_n \in \Lambda_N$ and either $\rank(a_1, \dots, a_n) < n$ or $\rank(a_1, \dots, a_n) = n$ and for any $1 \leq j \leq n$, $\alpha_j \mod L'$ is a generator of the cyclic group $L/L'$ (see Definition \ref{defgoodsmoothing}). Our main result is the following:

\begin{general}{Theorem}\label{maintheorem}
Suppose that $a_1, \dots, a_n \in \Lambda$ are non-zero linear forms which are \wellplaced in $\dualsimple{V}$ and assume that the smoothing lattice $L'$ is \good for $a_1, \dots, a_n$ where $n, N \geq 2$. Then there is an integer $b = b(a_1, \dots, a_n, v) \in \zz$ which depends on the linear forms $a_1, \dots, a_n$ and on the class of $v$ in $V/L'$ but not on $w, x \in \cc \times \homlc$ such that for all $(v, w, x) \in V/L' \times \cc \times \homlc$:
\begin{equation}\label{analogueDR}
\prod_{j = 1}^{n} G_{n-2, a_1, \dots, \omitvar{a_j}, \dots, a_n}(v)(w, \plgt, L, L')^{(-1)^{j+1}}= \exp\left(\frac{2i\pi b}{\denomclasscohom}\right)
\end{equation}
where $\denomclasscohom = \prod_{p | N} p^{\left\lfloor\frac{n}{p-1} \right\rfloor}$.
\end{general}

Formula \refp{analogueDR} is an analogue for general $n \geq 2$ of formula \refp{dedekindrademachersmooth}, which we may recover as follows. Fix $L = \zz e_1 + \zz e_2$ and $x \in \homlc$ satisfying $x(e_1) = \tau$ and $x(e_2) = 1$. Write $a_1 = [1,0]$ for the linear form on $L$ satisfying $a_1(e_1) = 1$ and $a_1(e_2) = 0$, as well as $a_2 = [d, -c]$ for the linear form satisfying $a_2(e_1) = d$ and $a_2(e_2) = -c$. Then, identifying $\theta = G_0$ we may write in our notations for $v_0 = 0$:
$$\theta_{[1, 0]}(v_0)(z, x) = \theta(z, \tau) ~~ \text{ and } ~~ \theta_{[d, -c]}(v_0)(z, x) = \theta\left(\frac{z}{c\tau + d}, \frac{a\tau + b}{c\tau + d}\right)$$
and the integer $b(a_1, a_2, v_0) = b([1,0], [d, -c], 0)$ is given by the formula:
$$P_{2, \gamma}^{(N)} = \frac{-b(a_1, a_2, v_0)}{\generaldenomclasscohom{N, 2}}$$
where $\gamma = \begin{psmallmatrix} \ast & \ast \\ c & d\end{psmallmatrix} \in \Gamma_0(N)$. In arithmetic applications, it is often useful to restrict the possible values of $z$ to the field $\qq(\tau)$ when $\tau$ is an real quadratic or imaginary quadratic number, in which case the regime $w = 0$, $z = v_1\tau + v_2$ with $v_1, v_2 \in \qq$ may be used (see for instance the determination of values of $L$ functions at $s=1$ for real quadratic fields in [\hspace{1sp}\cite{Siegel}, Chapter II, \symbolparagraph 6]). 

From Theorem \ref{maintheorem} we immediately deduce partial multiplicative cocycle relations for the smoothed $\smoothedgr$ raised to the power $\denomclasscohom$ as:

\begin{mycor}\label{cormaintheorem}
Under the same hypothesis as in Theorem \ref{maintheorem}:
$$\left(\prod_{j = 1}^{n} G_{n-2, a_1, \dots, \omitvar{a_j}, \dots, a_n}(v)(w, \plgt, L, L')^{(-1)^{j+1}}\right)^{\denomclasscohom} = 1.$$
\end{mycor}

We note that when the dimension $n$ is fixed, the integers $\denomclasscohom$ are also uniformly bounded by the integer $\mathcal{D}(n)$ defined by:
$$ \mathcal{D}(n) = \prod_{p \leq n+1}p^{\left\lfloor\frac{n}{p-1} \right\rfloor}.$$
For instance, when $n = 2$ we recover $\mathcal{D}(2) = 12$ so that for all $N \geq 2$ and for all $\gamma \in \Gamma_0(N)$, $\mathcal{D}(2).P_{2, \gamma}^{(N)} \in \zz$. The general bound $\mathcal{D}(n)$ appears in the study of certain higher dimensional Dedekind sums related to the smoothed functions $\smoothedbn$ (see section \ref{sectiontraces}) and it is a classical bound in the study of such objects (see \cite{ZagierDedekind} where other versions of higher dimensional Dedekind sums have denominators uniformly bounded for even $n$ by $\mu_{n/2} = 2^{-n}\mathcal{D}(n)$ in our notation).

Let us now give an outline of the paper, which is divided into two main parts. In section \ref{sectionmodularproperty} we recall the important definitions of \cite{firstpaper} and prove that relation \refp{modularproperty} holds for almost all configurations of the linear forms $a_1, \dots, a_n$, thus expanding [\hspace{1sp}\cite{firstpaper}, Theorem 1]:

\begin{general}{Theorem}\label{theoremmodularextend}
Let $a_1, \dots, a_n$ be non-zero linear forms in $\Lambda$ which are \wellplaced in $\dualsimple{V}$. Then:
$$\prod_{j = 1}^{n} G_{n-2, a_1, \dots, \omitvar{a_j}, \dots, a_n}(v)(w, \plgt)^{(-1)^{j+1}} = \exp(2i\pi B_{n, a_1, \dots, a_n}(v)(w,\plgt)).$$
\end{general}

In \cite{firstpaper} we already proved the case where $\rank(a_1, \dots, a_n) = n$ and the case where $\rank(a_1, \dots, a_n) \leq n-2$ is trivial. Thus, in this paper, we only prove the case where $\rank(a_1, \dots, a_n) = n-1$ and $0$ is not a barycenter of $a_1, \dots, a_n$ (see Proposition \ref{propmodularextend}). Theorem \ref{theoremmodularextend} is a first step in the direction of Theorem \ref{maintheorem}.

Next, in section \ref{sectionsmoothing} we apply to the $B_{n, a_1, \dots, a_n}$ functions a standard smoothing operation inspired by \cite{CD} and prove the following theorem:

\begin{general}{Theorem}\label{theoremsmoothing}
Let $a_1, \dots, a_n \in \Lambda$ be linearly independent and suppose that the smoothing lattice $L'$ is \good for $a_1, \dots, a_n$. Let $\alpha_1, \dots, \alpha_n$ be the primitive positive dual basis to $a_1, \dots, a_n$ in $L$. Let $\alpha_{1, j} = \langle \alpha_j, e_1 \rangle$, $s_j = a_j(\alpha_j)$ and $v = \sum_{j = 1}^n v_j\alpha_j/s_j$. Fix any set of representatives $\mathcal{F}$ for $L/M$ where $M = \oplus_{j = 1}^n \zz \alpha_j$. Then:
\begin{equation}\label{bernoullitraces}
\smoothedbn = \epsilon\sum_{\delta \in \mathcal{F}} \sum_{d | N, d \neq 1}\trace_{\qq(\zeta_d)/\qq}\left(\prod_{j = 1}^{n} \left(\frac{\zeta_d^{-\alpha_{1,j}\lfloor \frac{v_j + \delta_j}{s_j} \rfloor}}{\zeta_d^{\alpha_{1,j}}-1}\right)\right)
\end{equation}
where $\epsilon = \signdet(a_1,\dots, a_n)$, $\zeta_d = \exp(2i\pi/d)$ and $\delta = \sum_{j = 1}^n \delta_j \alpha_j/s_j$ for any $\delta \in \mathcal{F}$.
\end{general}

This theorem states that for fixed $a_1, \dots, a_n \in \Lambda$ and fixed $v \in V/L'$ the function $(w,x) \to  \smoothedbn$ is actually constant, provided that the smoothing lattice $L'$ is \good for $a_1, \dots, a_n$. In addition, for any $v \in V/L'$, the rational number $\smoothedbn$ is expressed as a sum of traces of algebraic numbers in cyclotomic fields whose denominators are well understood (see \cite{DasguptaShintani} or \cite{ZagierDedekind} for instance). It follows from a detailed analysis of these algebraic numbers:
$$\smoothedbn \in \denomclasscohom^{-1} \zz.$$
Theorem \ref{maintheorem} shall then be obtained as a consequence of this result together with Theorem \ref{theoremmodularextend}. 

This paper is organised as follows: in section \ref{sectionmodularproperty} we recall the main definitions of the first paper in this series and prove Theorem \ref{theoremmodularextend} by a careful analysis of the cones involved in the definition of the smoothed $\smoothedgr$ functions. In section \ref{sectionsmoothing} we perform the smoothing operation on the Bernoulli rational functions and prove Theorem \ref{theoremsmoothing}. Then, at the end of section \ref{sectionsmoothing} we prove Theorem \ref{maintheorem} as a consequence of Theorem \ref{theoremmodularextend} together with Theorem \ref{theoremsmoothing}. Lastly, in section \ref{sectionsmoothingcocycle} we give cohomological interpretations of our main results.

\section{The modular property}\label{sectionmodularproperty}

In this section we recall the definitions given in \cite{firstpaper} and define the general geometric setup for the rest of this paper (see section \ref{sectionrecall}) and then we give a five-step proof of Theorem \ref{theoremmodularextend} in section \ref{sectionproofmodularextend}.

\subsection{Geometric families of $G_{n-2}$ and $B_n$ functions}\label{sectionrecall}

In this section we fix the geometric setup and recall the definitions given in \cite{firstpaper}. Let $V$ be a $\qq$-vector space of finite dimension $n$ and $L$ be a rank $n$ lattice in $V$. We fix $B = [e_1, \dots, e_n]$ a $\zz$-basis of $L$. Define $\Lambda = \homlz$ and fix $B_{\Lambda} = [f_1, \dots, f_n]$ the basis of $\Lambda$ satisfying $f_j(e_k) = \delta_{jk}$ where $\delta_{jk}$ is Kronecker's symbol. This fixes determinant forms on both $L$ and $\Lambda$, as well as actions of $\slnz{n}$ on $L$ by left multiplication and on $\Lambda$ by inverse right multiplication such that $\forall\, (a, \alpha) \in \Lambda \times L$, $\forall\, g \in \slnz{n}$, $(g\cdot a)(g \cdot \alpha) = a(\alpha)$.
  
Let us now recall the definition of the collection of geometric $G_{n-2, a_1, \dots, a_{n-1}}$ functions attached to $n-1$ primitive linear forms in $\Lambda$. Assume that $a_1, \dots, a_{n-1}$ are linearly independent and fix $\alpha_1, \dots, \alpha_{n-1}$ a positive dual family to $a_1, \dots, a_{n-1}$ in $L$ in the sense of [\hspace{1sp}\cite{firstpaper}, Lemma 6], which means that for all $1 \leq j \leq n-1$:
$$a_k(\alpha_j) = 0 \text{ for all } k \neq j \text{ and } a_j(\alpha_j) > 0. $$
We may then define the set 
\begin{equation}\label{deffabaralphabarv}
F(\abar, \alphabar, v) = \{\delta \in v+L \setseparator 0 \leq a_j(\delta) < a_j(\alpha_j), \forall\, 1 \leq j \leq n-1 \}.
\end{equation}
as well as the primitive vector $\gamma \in L$ such that $\det(a_1, \dots, a_{n-1}, \cdot) = s\gamma$ for some integer $s>0$. It is clear that the set $F(\abar, \alphabar, v)$ is invariant under translation along $\gamma$ and it follows from [\hspace{1sp}\cite{firstpaper}, Proposition 7] that the function
\begin{equation}\label{originaldefinition}
G_{n-2, \abar}(v)(w,x) := \prod_{\delta \in F(\abar, \alphabar, v)/\zz\gamma}G_{n-2}\left(\frac{w + x(\delta)}{x(\gamma)}, \frac{x(\alpha_1)}{x(\gamma)}, \dots, \frac{x(\alpha_{n-1})}{x(\gamma)}\right)
\end{equation}
is well defined for $v \in V/L$ and $(w,x)$ in a dense open subset of $\cc \times \homlc \simeq \cc \times \cc^n$ and is independent of the choice of $\alpha_1, \dots, \alpha_{n-1}$. In the proof of [\hspace{1sp}\cite{firstpaper}, Proposition 7] we give an alternative definition of these functions which will be useful for the proof of Theorem \ref{theoremmodularextend} and which we now recall. Fix $x \in \homlc$ such that $x(\alpha_j)/x(\gamma) \not\in \rr$ for all $1 \leq j \leq n-1$. Then define the signs $d_j = \sign(\Im(x(\alpha_j)/x(\gamma))) \in \{-1, +1\}$ and set $D = \sum_{j = 1}^{n-1} (d_j-1)/2$. Finally, define the two cones:
\begin{align*}
C^{+}(\abar, x) &= \{\delta \in V \setseparator \forall\, 1 \leq j \leq n-1, a_j(\delta) \geq 0 \text{ if } d_j = 1, a_j(\delta) < 0 \text{ if } d_j = -1 \}, \\
C^{-}(\abar, x) &= \{\delta \in V \setseparator \forall\, 1 \leq j \leq n-1, a_j(\delta) \geq 0 \text{ if } d_j = -1, a_j(\delta) < 0 \text{ if } d_j = 1 \}.
\end{align*}
These two cones are independent from the choice of $\alpha_1, \dots, \alpha_{n-1}$ and are invariant by translation along $\gamma$. It follows from the proof of [\hspace{1sp}\cite{firstpaper}, Proposition 7] that:
\begin{multline}\label{alternativedefinition}
G_{n-2, a_1, \dots, a_{n-1}}(v)(w,x)^{(-1)^D} = \prod_{\delta \in (v+L) \cap C^{-}(\abar, x)/\zz\gamma} \left(1-e^{-2i\pi\left(\frac{w+\plgt(\delta)}{\plgt(\gamma)}\right)}\right)\\ \times \prod_{\delta \in (v+L) \cap C^{+}(\abar, x)/\zz\gamma} \left(1-e^{2i\pi\left(\frac{w+\plgt(\delta)}{\plgt(\gamma)}\right)}\right)^{(-1)^n}.
\end{multline}
This second formulation will be used for the proof of Theorem \ref{theoremmodularextend} in section \ref{sectionproofmodularextend}. Finally, when the linear forms $a_1, \dots, a_{n-1}$ are not linearly independent, we define by convention $G_{n-2, a_1, \dots, a_{n-1}}(v)(w,x) = 1$ for any $v, w, x \in V/L \times \cc \times \homlc$. 

Let us now recall the definition of the associated family of $B_{n, a_1, \dots, a_n}$ functions. If the $n$ primitive linear forms $a_1, \dots, a_n \in \Lambda$ are linearly dependent, we set $B_{n, a_1, \dots, a_n}(v)(w,x) = 0$. Otherwise, there is a unique family $\alpha_1, \dots, \alpha_n$ of primitive vectors in $L$ satisfying for all $1 \leq j \leq n$:
$$a_k(\alpha_j) = 0 \text{ for all } k \neq j \text{ and } a_j(\alpha_j) > 0. $$
Set $P(\abar) = \{ \delta \in V \setseparator \delta = \sum_{j = 1}^n \delta_j \alpha_j, 0 \leq \delta_j < 1 \text{ for all } 1 \leq j \leq n\}$. Then the function $B_{n, a_1, \dots, a_n} : V/L \to \qq[w](x)$ is defined by:
\begin{equation}\label{defbernoulliexplicit}
B_{n, a_1, \dots, a_n}(v)(w,x) = \signdet(a_1, \dots, a_n)\times \coefficient{\sum_{\delta \in (v+L)\cap P(\abar)}\frac{e^{wt}e^{x(\delta)t}}{\prod_{j = 1}^n (1 - e^{x(\alpha_j)t})}}{t}{0}.
\end{equation}
In [\hspace{1sp}\cite{firstpaper}, Proposition 7] we proved that if $a_1, \dots, a_n$ are linearly independent, then:
$$ \prod_{j= 1}^{n} G_{n-2, a_1, \dots, \omitvar{a_j}, \dots, a_n}(v)(w, x)^{(-1)^{j+1}} = \exp(2i\pi B_{n, a_1, \dots, a_n}(v)(w,x))$$
In addition, the two collections of functions $G_{n-2, a_1, \dots, a_{n-1}}$ and $B_{n, a_1, \dots, a_n}$ are equivariant under the action of $\slnz{n}$ on $L$ as for all $g \in \slnz{n}$:
\begin{align*}
G_{n-2,g \cdot a_1, \dots, g \cdot a_{n-1}}(g \cdot v)(w, g\cdot x) &= G_{n-2, a_1, \dots, a_{n-1}}(v)(w, x), \\
B_{n, g\cdot a_1, \dots, g\cdot a_n}(g \cdot v)(w, g\cdot x) &= B_{n, a_1, \dots, a_n}(v)(w, x).
\end{align*}

Before giving the proof of Theorem \ref{theoremmodularextend} we recall the definition of the standard non-trivial relation given in [\hspace{1sp}\cite{firstpaper}, Definition 11] in regards to the condition on the positions of the linear forms $a_1, \dots, a_n$ under which we formulate Theorem \ref{theoremmodularextend}. Indeed, if $a_1, \dots, a_m$ are $m$ non-zero linear forms on $V$ such that $\rank(a_1, \dots, a_m) = m-1$ then the standard non-trivial relation among $a_1, \dots, a_m$ is the unique relation $\sum_{j = 1}^m \lambda_j a_j = 0$ with rational coefficient $\lambda_1, \dots, \lambda_n$ satisfying:
\begin{itemize}
\item there is an index $1 \leq l \leq n$ such that $\lambda_j = 0$ for $1 \leq j < l$ and $|\lambda_l| = 1$.
\item $\cardinalshort{\{1 \leq j \leq n \setseparator \lambda_j < 0 \}} \leq \cardinalshort{\{1 \leq j \leq n \setseparator \lambda_j > 0\}}.$
\item in case of equality in the previous condition, $\lambda_l = -1$.
\end{itemize}
In this situation we define, as in \cite{firstpaper}, $k^{\pm}(a_1, \dots, a_m) =  \cardinalshort{\{ 1 \leq j \leq m \setseparator \pm \lambda_j > 0\}}$. The condition on the positions of the linear forms $a_1, \dots, a_n$ we briefly presented in the introduction may be stated in terms of this standard non-trivial relation:

\begin{general}{Definition}\label{definitiongoodposition}
Let $a_1, \dots, a_n$ be $n$ linear forms on $V$. We say that $a_1, \dots, a_n$ are \wellplaced in $\dualsimple{V}$ if either $\rank(a_1, \dots, a_n) \neq  n-1$ or if $\rank(a_1, \dots, a_n) = n-1$ and $0$ is not a barycenter of $a_1, \dots, a_n$. Equivalently, $a_1, \dots, a_n$ are \wellplaced in $\dualsimple{V}$ if $\rank(a_1, \dots, a_n) \neq n-1$ or if $\rank(a_1, \dots, a_n) = n-1$ and $k^{-}(a_1, \dots, a_n) > 0$.
\end{general}

Notice that this condition is very similar to the ``good position'' condition in \cite{firstpaper} for $n+1$ linear forms under which the cocycle relation \refp{cocyclebernoulli} holds. The rest of section \ref{sectionmodularproperty} is devoted to the proof of Theorem \ref{theoremmodularextend} in the case where $\rank(a_1, \dots, a_n) = n-1$ and $k^{-}(a_1, \dots, a_n) > 0$ and will make use of the standard non-trivial relation among $a_1, \dots, a_n$ which we have recalled.

\subsection{Proof of the modular property}\label{sectionproofmodularextend}

In this section we give the proof of Theorem \ref{theoremmodularextend}. We once again highlight that when $\rank(a_1, \dots, a_n) = n$, the theorem was proven in \cite{firstpaper} using Narukawa's theorem (see \cite{Narukawa}), and that when $\rank(a_1, \dots, a_n) \leq n-2$ the statement is trivial. We therefore only need to prove the following.

\begin{general}{Proposition}\label{propmodularextend}
Let $a_1, \dots, a_n \in \Lambda$ be $n$ non-zero primitive linear forms such that $\rank(a_1, \dots, a_n) = n-1$ and $k^{-}(a_1, \dots, a_n)>0$. Then for all $v, w, x \in V/L \times \cc \times \homlc$:
\begin{equation}\label{formulamodularextend}
\prod_{j = 1}^n G_{n-2, a_1, \dots, \omitvar{a_j}, \dots, a_n}(v)(w,x)^{(-1)^{j+1}} = 1.
\end{equation}
\end{general}

\noindent \textbf{Remark:} If we remove the assumption that $k^{-}(a_1, \dots, a_n) > 0$ the result does not generally hold. To see this we analyse the simple case where $a_1 = -a_{2}$ and $\rank(a_2, \dots, a_{n}) = n-1$ with $n \geq 3$. In this case, the left-hand side of \refp{formulamodularextend} reduces to:
$$\prod_{j = 1}^n G_{n-2, a_1, \dots, \omitvar{a_j}, \dots, a_n}(v)(w,x)^{(-1)^{j+1}}  = \frac{G_{n-2, -a_2, a_3, \dots, a_n}(v)(w,x)}{G_{n-2, a_2, a_3, \dots, a_n}(v)(w,x)}.$$
For simplicity, fix a $\zz$-basis $B = [e_1, \dots, e_n]$ of $L$ and fix $a_j$ the linear form satisfying $a_j(e_k) = 0$ if $k \neq j$, $a_j(e_j) = 1$ for $2 \leq j \leq n$. If we assume further that $F(\abar, \alphabar, v) = \{v\}$ where $\abar = (a_2, \dots, a_n)$ and $\alphabar = (e_2, \dots, e_n)$ (see \refp{deffabaralphabarv} for the definition of $F(\abar, \alphabar, v)$) and that $a_2(v) = 0$ then:
$$\frac{G_{n-2, -a_2, a_3, \dots, a_n}(v)(w,x)}{G_{n-2, a_2, a_3, \dots, a_n}(v)(w,x)} = \frac{G_{n-2}\left(\frac{w + x(v)}{x(-e_1)}, \frac{x(-e_2)}{x(-e_1)}, \frac{x(e_3)}{x(-e_1)}, \dots, \frac{x(e_n)}{x(-e_1)}\right)}{G_{n-2}\left(\frac{w + x(v)}{x(e_1)}, \frac{x(e_2)}{x(e_1)}, \frac{x(e_3)}{x(e_1)}, \dots, \frac{x(e_n)}{x(e_1)}\right)}.$$
Using [\hspace{1sp}\cite{Nishizawa}, Proposition 3.2] we get:
$$ \frac{G_{n-2, -a_2, a_3, \dots, a_n}(v)(w,x)}{G_{n-2, a_2, a_3, \dots, a_n}(v)(w,x)} = G_{n-3}\left(\frac{w + x(v)}{x(e_1)}, \frac{x(e_3)}{x(e_1)}, \dots, \frac{x(e_n)}{x(e_1)}\right)$$
and this is not identically equal to $1$. For $n = 2$ we get under the same assumptions the simpler form
$$\theta_{-a}(v)(w,x)\theta_{a}(w,x)^{-1} = \exp\left(-2i\pi\left(\frac{w+x(v)}{x(e_1)} - \frac12\right)\right)$$
which is also not identically equal to $1$. \smallskip

\noindent We organise the proof of Proposition \ref{propmodularextend} into five main steps:
\begin{itemize}
\item \underline{Step 1:} we first show in section \ref{modextendpartone} that we can order the linear forms $a_1, \dots, a_n$ at will. To simplify the notations, we shall choose an order on $a_1, \dots, a_n$ such that the standard non-trivial relation among $a_1, \dots, a_n$ is $\sum_{j = 1}^n \lambda_j a_j = 0$ with 
$$ \begin{cases} \lambda_j = 0 & \text{ for } 1 \leq j < l \\ 
\lambda_j < 0 & \text{ for } l \leq j < m \\
\lambda_j > 0 & \text{ for } m \leq j \leq n 
\end{cases} $$
for some $1 \leq l < m \leq n$.
\item \underline{Step 2:} In section \ref{modextendparttwo}, we shall rewrite each term in the right-hand side of \refp{formulamodularextend} using formula \refp{alternativedefinition} and reorganise the factors. More precisely, we define two families of cones $(\Cun{j})_j$ and $(\Cdeux{j})_j$ for $l \leq j \leq n$ (see Definition \ref{defcuncdeux}) and prove in Lemma \ref{lemmafunfdeux} that: 
\begin{multline}\label{sketchproofcones}
\prod_{j = 1}^n G_{n-2, a_1, \dots, \omitvar{a_j}, \dots, a_n}^{(-1)^{j+1}}(v)(w,x) = \prod_{j = l}^n \prod_{\delta \in (v+L) \cap \Cun{j}/\zz\gamma} \left(1-e^{2i\pi\left(\frac{w+\plgt(\delta)}{\plgt(\gamma)}\right)}\right)^{\mu_j} \\ \times \prod_{\delta \in (v+L) \cap \Cdeux{j}/\zz\gamma} \left(1-e^{-2i\pi\left(\frac{w+\plgt(\delta)}{\plgt(\gamma)}\right)}\right)^{\mu_j (-1)^n}
\end{multline}
where $\gamma$ is the primitive vector in $L$ satisfying $\det(a_1, \dots, a_{n-1}, \cdot) = s_n\gamma$ for some positive integer $s_n$ and the $\mu_j$'s are explicit signs in $\{-1, +1\}$. Denoting by $\cun{j}$ and $\cdeux{j}$ the indicator functions associated to these cones we may define two functions $\fun, \fdeux : V/\qq\gamma \to \zz$ by:
$$\fun = \sum_{j = l}^n \mu_j \cun{j} ~~ \text{ and } ~~ \fdeux = \sum_{j = l}^n \mu_j \cdeux{j} \times (-1)^n.$$
Then formula \refp{sketchproofcones} may be rewritten as:
\begin{equation}\label{sketchprooffunfdeux}
\prod_{j = 1}^n G_{n-2, a_1, \dots, \omitvar{a_j}, \dots, a_n}(v)(w,x)^{(-1)^{j+1}} = \prod_{\delta \in v+L/\zz\gamma} \left(1-e^{2i\pi\left(\frac{w+\plgt(\delta)}{\plgt(\gamma)}\right)}\right)^{\fun(\delta)}\left(1-e^{-2i\pi\left(\frac{w+\plgt(\delta)}{\plgt(\gamma)}\right)}\right)^{\fdeux(\delta)}.
\end{equation}
\item \underline{Step 3:} The remainder of the proof is devoted to showing that $\fun = \fdeux = 0$. The proof that $\fdeux = 0$ is exactly the same as the proof that $\fun = 0$ on which we now focus. This is done by a technical combinatorial analysis of the cones $\Cun{j}$. In section \ref{modextendpartthree} we give an example for $n = 4$ where we show how the table of the signs of the linear forms $a_k$ on the cones $\Cun{j}$ contains the relevant information for the proof that $\fun = 0$. We shall show that this table of signs must obey certain rules (see Lemma \ref{lemmasignrelations}), for instance the signs of $a_k$ on $\Cun{j}$ must be related to the sign of $a_j$ on $\Cun{k}$ for $l \leq j \neq k \leq n$. We then deduce from this set of rules that if $l\leq j, j', j'' \leq n$ are three  distinct indices then the triple intersection $\Cun{j} \cap \Cun{j'} \cap \Cun{j''}$ is empty (see Lemma \ref{lemmatroisvide}).
\item \underline{Step 4:} In section \ref{modextendpartfour} we show that any vector $\delta \in V$ belongs to exactly $0$ or $2$ of the cones $\Cun{l}, \dots, \Cun{n}$. The third step in the proof guarantees that any $\delta \in V$ belongs to either $0$, $1$ or $2$ of these cones, so we only need to show (see Lemma \ref{lemmastepfive}) that a vector $\delta \in V$ cannot belong to exactly one of these cones. This is by far the most technical part of the proof, relying on a technical property of the sign table (see Lemma \ref{lemmasignrelations}, (iii)).
\item \underline{Step 5:} In section \ref{modextendpartfive} we complete the proof by showing that if $\delta \in \Cun{j} \cap \Cun{j'}$ for some $j \neq j'$ then $\fun(\delta) = 0$. This is done by showing that in this case $\mu_j = -\mu_{j'}$ (see Lemma \ref{lemmafunzero}). 
\end{itemize}

\subsubsection{Invariance under permutation}\label{modextendpartone}

In this section we justify that Proposition \ref{propmodularextend} holds for the linear forms $a_1, \dots, a_n$ if and only if it holds for any permutation $a_{\sigma(1)}, \dots, a_{\sigma(n)}$ of $a_1, \dots, a_n$ where $\sigma \in \goth{S}_n$ using the following result.

\begin{general}{Lemma}
Let $a_1, \dots, a_n \in \Lambda$ be $n$ non-zero linear forms. For any permutation $\sigma \in \goth{S}_{n}$:
$$\prod_{j = 1}^{n} G_{n-2, a_{\sigma(1)}, \dots, \omitvar{a_{\sigma(j)}}, \dots, a_{\sigma(n)}}^{(-1)^{j+1}} = \left(\prod_{j = 1}^n G_{n-2, a_1, \dots, \omitvar{a_j}, \dots, a_n}^{(-1)^{j+1}}\right)^{\signature(\sigma)}$$
\end{general}

\begin{proof}
As the transpositions generate $\goth{S}_n$ it is sufficient to prove this statement for transpositions. Fix $\sigma = (kl)$ the transposition switching $k, l$ with $k < l$. We wish to prove that:
$$\prod_{j = 1}^{n} G_{n-2, a_{\sigma(1)}, \dots, \omitvar{a_{\sigma(j)}}, \dots, a_{\sigma(n)}}^{(-1)^{j+1}} = \left(\prod_{j = 1}^n G_{n-2, a_1, \dots, \omitvar{a_j}, \dots, a_n}^{(-1)^{j+1}}\right)^{-1}$$
We will repeatedly use the fact that if $b_1, \dots, b_{n-1}$ are non-zero linear forms then for any permutation $\rho \in \goth{S}_{n-1}$:
\begin{equation}\label{permutationgr}
G_{n-2, b_{\rho(1)}, \dots, b_{\rho(n-1)}} = G_{n-2, b_1, \dots, \omitvar{b_j}, \dots, b_n}^{\signature(\rho)}
\end{equation}
which is clear from the definition (see \refp{originaldefinition}). Consider first an index $1 \leq j \leq n$ such that $j \neq k$ and $j \neq l$. Then $\sigma$ reduces to the transposition $(kl)$ on the set $\{1, \dots, j-1, j+1, \dots, n\}$ and it follows from \refp{permutationgr} that:
$$G_{n-2, a_{\sigma(1)}, \dots, \omitvar{a_{\sigma(j)}}, \dots, a_{\sigma(n)}} = G_{n-2, a_1, \dots, \omitvar{a_j}, \dots, a_n}^{-1}.$$
Suppose now that $j = k$.
Then $$G_{n-2, a_{\sigma(1)}, \dots, \omitvar{a_{\sigma(k)}}, \dots, a_{\sigma(n)}} = G_{n-2, a_1, \dots, a_{k-1}, a_{k+1}, \dots, a_{l-1}, a_k, a_{l+1}, \dots, a_n} = G_{n-2, a_{\rho(1)}, \dots, \omitvar{a_{\rho(l)}}, \dots, a_{\rho(n)}}$$
where $\rho$ is the cycle $(k, l, l-1, \dots, k+1)$ which has signature $(-1)^{l+k+1}$. Therefore, formula \refp{permutationgr} implies that
$$ G_{n-2, a_{\sigma(1)}, \dots, \omitvar{a_{\sigma(k)}}, \dots, a_{\sigma(n)}}^{(-1)^{k+1}} = G_{n-2, a_1, \dots, \omitvar{a_l}, \dots, a_n}^{(-1)^{l}}$$
and we may prove similarly if $j = l$ that:
$$ G_{n-2, a_{\sigma(1)}, \dots, \omitvar{a_{\sigma(l)}}, \dots, a_{\sigma(n)}}^{(-1)^{l+1}} = G_{n-2, a_1, \dots, \omitvar{a_k}, \dots, a_n}^{(-1)^{k}}$$
This gives the desired result:
$$ \prod_{j = 1}^{n} G_{n-2, a_{\sigma(1)}, \dots, \omitvar{a_{\sigma(j)}}, \dots, a_{\sigma(n)}}^{(-1)^{j+1}} = \left(\prod_{j = 1}^n G_{n-2, a_1, \dots, \omitvar{a_j}, \dots, a_n}^{(-1)^{j+1}}\right)^{-1}.$$ 
\end{proof}

In particular to show Proposition \ref{propmodularextend} we may switch the order of the linear forms $a_1, \dots, a_n$. Thus, for the remainder of section \ref{sectionproofmodularextend} we assume that $a_1, \dots, a_n$ are non-zero linear forms such that $\rank(a_1, \dots, a_n) = n-1$ and such that the standard non-trivial relation $\sum_{j = 1}^n \lambda_j a_j = 0$ among $a_1, \dots, a_n$ (see [\hspace{1sp}\cite{firstpaper}, Definition 11]) satisfies:
\begin{equation}\label{signslambda}
\begin{cases} \lambda_j = 0 & \text{ for } 1 \leq j < l \\ 
\lambda_j < 0 & \text{ for } l \leq j < m \\
\lambda_j > 0 & \text{ for } m \leq j \leq n 
\end{cases} 
\end{equation}
for some $1 \leq l < m \leq n$. We end this section by making the following remark: there are two simpler cases corresponding to $l = n-1$ and to $l = 1$. If $l = n-1$ then the right-hand side of \refp{formulamodularextend} is 
$$ G_{n-2, a_1, \dots, a_{n-1}}(v)(w,x)^{(-1)^{n}}G_{n-2, a_1, \dots, a_{n-1}}(v)(w,x)^{(-1)^{n+1}} = 1. $$
In the case where $l = 1$, formula \refp{formulamodularextend} may be obtained by introducing a linear form $a_{n+1}$ such that $\rank(a_1, \dots, a_{n+1}) = n$ and by applying directly [\hspace{1sp}\cite{firstpaper}, Theorem 1] in conjunction with [\hspace{1sp}\cite{firstpaper}, Corollary 2] to the families $a_1, \dots, \omitvar{a_j}, \dots, a_{n+1}$ for $1 \leq j \leq n$. This strategy however doesn't generalise to the case $1 < l < n-1$ for which we need the proof presented in this paper.

\subsubsection{Definition of the cones $\Cun{j}$ and $\Cdeux{j}$}\label{modextendparttwo}

We now explicitly describe the product in the left-hand side of \refp{formulamodularextend} in the case where the coefficients $\lambda_j$ satisfying $\sum_{j = 1}^n \lambda_j a_j = 0$ verify \refp{signslambda}. We remark that whenever $1 \leq j < l$ we get $\rank(a_1, \dots, \omitvar{a_j}, \dots, a_n) = n-2$ therefore $ G_{n-2, a_1, \dots, \omitvar{a_j}, \dots, a_n} = 1$ by definition. In that case, formula \refp{formulamodularextend} which we aim to prove reduces to:
$$\prod_{j = l}^{n} G_{n-2, a_1, \dots, \omitvar{a_j}, \dots, a_n}(v)(w, \plgt)^{(-1)^{j+1}} = 1$$ 
Let us now describe each of the non-trivial terms in the left-hand side of \refp{formulamodularextend}, using the alternative definition \refp{alternativedefinition} for the function $G_{n-2, a_1, \dots, \omitvar{a_j}, \dots, a_n}$. Indeed, for any fixed $l \leq j \leq n$ define 
$\gammaj{j}$ to be the unique primitive vector in $L$ such that $\det(a_1, \dots, \omitvar{a_j}, \dots, a_n, \cdot) = s^{(j)}\gammaj{j}$ for some positive integer $s^{(j)}$. Fix $(\alphajk{j}{k})_{k \neq j}$ a positive dual family to $a_1, \dots, \omitvar{a_j}, \dots, a_n$ in $L$, i.e. a family satisfying for all $1 \leq k \leq n$, $k \neq j$:
$$ a_k(\alphajk{j}{k'}) = 0 \text{ for all } 1 \leq k' \leq n, k' \neq k, j \text{ and } a_k(\alphajk{j}{k}) > 0.$$
Next we fix $\plgt \in \homlc$ such that for all $l \leq j \leq n$ and all $1 \leq k \leq n$, $k \neq j$, we have $\plgt(\alphajk{j}{k})/\plgt(\gammaj{j}) \not\in \rr$. Let us then define the linear form $\yj{j} : V \to \rr$ by $\yj{j}(v) = \Im(x(v)/x(\gammaj{j}))$. For $1 \leq k \leq n$, $k \neq j$, define 
\begin{equation}\label{defdjk}
\djk{j}{k} = \sign(\yj{j}(\alphajk{j}{k})) \in \{\pm1\}
\end{equation}
and set $D_j = \sum_{k \neq j} (\djk{j}{k}-1)/2$. Finally, recall from section \ref{sectionrecall} the cones:
\begin{equation}\label{defcplusmoins}
C^{\pm}_j = \{\delta \in V \setseparator \forall\, 1 \leq k \leq n, k \neq j,  a_k(\delta) \geq 0 \text{ if } \pm\djk{j}{k} = 1, a_k(\delta) < 0 \text{ if } \pm\djk{j}{k} = -1\}.
\end{equation}
It follows from formula \refp{alternativedefinition} that for any $l \leq j \leq n$:
\begin{multline}\label{eqgrdefinition}
G_{n-2, a_1, \dots, \omitvar{a_j}, \dots, a_n}(v)(w, \plgt) = \prod_{\delta \in (v+L) \cap C^{-}_j/\zz\gammaj{j}} \left(1-e^{-2i\pi\left(\frac{w+\plgt(\delta)}{\plgt(\gammaj{j})}\right)}\right)^{(-1)^{D_j}}\\ \times \prod_{\delta \in (v+L) \cap C^{+}_j/\zz\gammaj{j}} \left(1-e^{2i\pi\left(\frac{w+\plgt(\delta)}{\plgt(\gammaj{j})}\right)}\right)^{(-1)^{D_j + n}}
\end{multline}
Let us now fix $\gamma = \gammaj{n}$, so that for any $l \leq j \leq n$, $\gammaj{j} = (-1)^{j+n}\sign(\lambda_j)\gamma$, where we recall that $\sum_{j = l}^n \lambda_j a_j = 0$ is the standard non-trivial relation among $a_1, \dots, a_n$. For simplicity, we shall define the sign
\begin{equation}\label{defepsj}
\eps_j = (-1)^{j+n}\sign(\lambda_j) \in \{-1, +1\}.
 \end{equation} 
 The reorganisation of the terms in \refp{eqgrdefinition} will be made by relabeling the cones $C^{+}_j$ and $C^{-}_j$ depending on the value of $\eps_j$:

\begin{general}{Definition}\label{defcuncdeux}
For $l \leq j \leq n$, define:
$$(\Cun{j}, \Cdeux{j}) = \begin{cases} (C^{+}_j, C^{-}_j) & \text{ if } \eps_j = 1 \\ (C^{-}_j, C^{+}_j) & \text{ if } \eps_j = -1 \end{cases} $$
and denote by $\cun{j}, \cdeux{j} : V \to \{0, 1\}$ their indicator functions. 
\end{general}

Note that we have explicitly:
\begin{align*}
\Cun{j} &= \{ \delta \in V \setseparator \forall 1 \leq k \leq n, k \neq j, a_k(\delta) \geq 0\text{ if } \eps_j\djk{j}{k} = 1, a_k(\delta) < 0 \text{ if } \eps_j\djk{j}{k} = -1\} \\
\Cdeux{j} &= \{ \delta \in V \setseparator \forall 1 \leq k \leq n, k \neq j, a_k(\delta) \geq 0\text{ if } \eps_j\djk{j}{k} = -1, a_k(\delta) < 0 \text{ if } \eps_j\djk{j}{k} = 1\}.
\end{align*}
Now, since $a_k(\gamma) = 0$ for any $1 \leq k \leq n$, it is clear that $\cun{j}(v + m\gamma) = \cun{j}(v)$ and $\cdeux{j}(v+m\gamma) = \cdeux{j}(v)$ for any $v \in V$ and any $m \in \qq$. Therefore, both functions $\cun{j}, \cdeux{j}$ reduce to functions on the quotient space $V/\qq\gamma$. Finally we define the signs:
\begin{equation}\label{defmuj}
\mu_j = \begin{cases} (-1)^{j+1+D_j+n} & \text{ if } \eps_j = 1 \\ (-1)^{j+1+D_j} & \text{ if } \eps_j = -1\end{cases}
\end{equation}
We are now ready to give a simple form for the left-hand side of formula \refp{formulamodularextend} which we will use for the rest of the proof.

\begin{general}{Lemma}\label{lemmafunfdeux}
With notations as above:
$$\prod_{j = l}^{n} G_{n-2, a_1, \dots, \omitvar{a_j}, \dots, a_n}(v)(w, \plgt)^{(-1)^{j+1}} = \prod_{\delta \in (v+L)/\zz\gamma}\left(1-e^{2i\pi\left(\frac{w+\plgt(\delta)}{\plgt(\gamma)}\right)}\right)^{\fun(\delta)} \left(1-e^{-2i\pi\left(\frac{w+\plgt(\delta)}{\plgt(\gamma)}\right)}\right)^{\fdeux(\delta)}$$
where the functions $\fun, \fdeux : V/\qq\gamma \to \zz$ are defined by:
$$\fun = \sum_{j = l}^n \mu_j \cun{j} ~~ \text{ and } ~~ \fdeux = \sum_{j = l}^n \mu_j \cdeux{j} \times (-1)^n.$$
\end{general}

\begin{proof}
Let us briefly denote by $J^{\pm}$ the set of indices $l \leq j \leq n$ satisfying $\eps_j = \pm 1$. Let us rewrite the terms in formula \refp{eqgrdefinition} for each $j \in J^+$ as:
\begin{multline*}
\prod_{j \in J^{+}} G_{n-2, a_1, \dots, \omitvar{a_j}, \dots, a_n}(v)(w, \plgt)^{(-1)^{j+1}} = \prod_{j \in J^{+}} \prod_{\delta \in (v+L) \cap C^{-}_j/\zz\gamma} \left(1-e^{-2i\pi\left(\frac{w+\plgt(\delta)}{\plgt(\gamma)}\right)}\right)^{\mu_j(-1)^n}\\ \times \prod_{\delta \in (v+L) \cap C^{+}_j/\zz\gamma} \left(1-e^{2i\pi\left(\frac{w+\plgt(\delta)}{\plgt(\gamma)}\right)}\right)^{\mu_j}.
\end{multline*}
On the other hand, since for all $j \in J^-$, $\gammaj{j} = -\gamma$, the product over $j \in J^-$ is:
\begin{multline*}
\prod_{j \in J^{-}} G_{n-2, a_1, \dots, \omitvar{a_j}, \dots, a_n}(v)(w, \plgt)^{(-1)^{j+1}} = \prod_{j \in J^{-}} \prod_{\delta \in (v+L) \cap C^{-}_j/\zz\gamma} \left(1-e^{2i\pi\left(\frac{w+\plgt(\delta)}{\plgt(\gamma)}\right)}\right)^{\mu_j}\\ \times \prod_{\delta \in (v+L) \cap C^{+}_j/\zz\gamma} \left(1-e^{-2i\pi\left(\frac{w+\plgt(\delta)}{\plgt(\gamma)}\right)}\right)^{\mu_j(-1)^n}.
\end{multline*}
Putting everything together and using the relabeled cones $\Cun{j}$ and $\Cdeux{j}$ gives:
\begin{multline*}
\prod_{j = l}^n G_{n-2, a_1, \dots, \omitvar{a_j}, \dots, a_n}(v)(w, \plgt)^{(-1)^{j+1}} = \prod_{j = l}^n \prod_{\delta \in (v+L) \cap \Cun{j}/\zz\gamma} \left(1-e^{2i\pi\left(\frac{w+\plgt(\delta)}{\plgt(\gamma)}\right)}\right)^{\mu_j}\\ \times \prod_{\delta \in (v+L) \cap \Cdeux{j}/\zz\gamma} \left(1-e^{-2i\pi\left(\frac{w+\plgt(\delta)}{\plgt(\gamma)}\right)}\right)^{\mu_j(-1)^n}
\end{multline*}
Thus the functions $\fun$ and $\fdeux$ are defined precisely so that:
$$\prod_{j = l}^{n} G_{n-2, a_1, \dots, \omitvar{a_j}, \dots, a_n}(v)(w, \plgt)^{(-1)^{j+1}} = \prod_{\delta \in (v+L)/\zz\gamma}\left(1-e^{2i\pi\left(\frac{w+\plgt(\delta)}{\plgt(\gamma)}\right)}\right)^{\fun(\delta)} \left(1-e^{-2i\pi\left(\frac{w+\plgt(\delta)}{\plgt(\gamma)}\right)}\right)^{\fdeux(\delta)}$$
which is the desired result.
\end{proof}

\subsubsection{Sign tables and emptyness of triple intersections}\label{modextendpartthree}

The rest of the proof of Proposition \ref{propmodularextend} consists in proving that $\fun = 0$ and $\fdeux = 0$. Both statements are proven similarly so we focus on the proof that $\fun = 0$.

Let us briefly explain the general idea on a simple example in the four-dimensional case. Let us fix:
$$a_1 = [1, 0, 0, 0], ~a_2 = [0, 1, 0, 0], ~a_3 = [0, 0, -1, 0], ~a_4 = [0, 1, 1, 0]. $$
These are four linear forms on $\zz^4$ such that $\rank(a_1, a_2, a_3, a_4) = 3$ and $0\cdot a_1 - a_2 + a_3 + a_4 = 0$ so that $l = 2$ and $m = 3$. The vector $\gamma$ is given by $\gamma = [0, 0, 0, -1]^T$. Let us also fix $x = [2i, 3i, 5i, -1]$ so that $x(\gamma) = 1$. We wish to describe explicitly the cones $C^1_j$ for $2 \leq j \leq 4$, therefore we need to describe the signs $\eps_j\djk{j}{k}$ for $2 \leq j \leq 4$ and $1 \leq k \leq 4$, $k \neq j$. We may compute a choice of elements $\alphajk{j}{k}$ as follows:

\newcommand\graycell{\cellcolor{gray}}

\begin{table}[H]
\begin{center}
\begin{tabular}{|C|C|C|C|}
\hline
$\diaghead{Very long}%
	{k}{j}$ 
	& 2 & 3 & 4\\
\hline
1 & \alphajk{2}{1} =[1, 0, 0, 0]^T & \alphajk{3}{1} = [1, 0, 0, 0]^T & \alphajk{4}{1} = [1, 0, 0, 0]^T \\ 
\hline 
2 & \graycell & \alphajk{3}{2} =  [0, 1, -1, 0]^T & \alphajk{4}{2} = [0, 1, 0, 0]^T \\ 
\hline
3 & \alphajk{2}{3} = [0, 1, -1, 0]^T  & \graycell & \alphajk{4}{3} = [0, 0, -1, 0]^T\\ 
\hline
4 & \alphajk{2}{4} = [0, 1, 0, 0]^T & \alphajk{3}{4} = [0, 0, 1, 0]^T  & \graycell \\ 
\hline
\end{tabular}
\end{center}
\caption{Example of a table of the $\alphajk{j}{k}$s}
\end{table}

\noindent This gives the computation of the signs $\eps_j\djk{j}{k}$:

\begin{table}[H]
\begin{center}
\begin{tabular}{|C|C|C|C|}
\hline
$\diaghead{Very long}%
  {k}{j}$
   & 2 & 3 & 4\\
\hline
1 & + &+ & +\\  
\hline
2 & \graycell & - & + \\ 
\hline
3 & -  & \graycell & -\\ 
\hline
4 &+ & + & \graycell \\ 
\hline
\end{tabular}
\end{center}
\caption{Example of a sign table containing the $\eps_j \djk{j}{k}$s}
\label{signtable}
\end{table}

From this sign table we deduce that $\Cun{3} \cap \Cun{4} = \emptyset$ as the conditions on $\eps_2\djk{2}{3} = - \eps_2\djk{2}{4}$ are incompatible. If we set $H^{+}_j = \{\delta \in V \setseparator a_j(\delta) \geq 0\}$ and $H^{-}_j = V - H^{+}_j$ then it follows from this table that:
$$\Cun{2} \cap \Cun{3} = \Cun{2} \cap H^{-}_2 = \Cun{3} \cap H^{-}_3 $$
and
$$\Cun{2} \cap \Cun{4} = \Cun{2} \cap H^{+}_2 = \Cun{4} \cap H^{+}_4 $$
In addition, the relation $a_2 = a_3 + a_4$ implies that $\Cun{3} \cap H^{+}_3 = \emptyset$ and $\Cun{4} \cap H^{-}_4 = \emptyset$. Therefore, $\Cun{2} = \Cun{3} \disjointunion \Cun{4}$ and we may check that $\fun = -\cun{2} + \cun{3} + \cun{4} = 0$. This is the general idea of the proof and we prove the general case in what follows.

As showcased by the previous example, we need to study the cones $\Cun{j}$ and therefore the signs $\eps_j\djk{j}{k}$ for $k \neq j$. We start by proving a crucial lemma on the relations between these signs that govern the sign tables (see Table \ref{signtable}). This will be useful for the last three steps of the proof of Proposition \ref{propmodularextend}. As a corollary, we will prove that any intersection of three of the $\Cun{j}$s is empty.

There are three main relations between the signs $\eps_j \djk{j}{k}$ among which the first two are quite simple. The third one is more technical and to state it we need to introduce for $l \leq j \leq n$ and $1 \leq k \leq n$, $k \neq j$ the following positive real numbers:
\begin{equation}\label{defujk}
\ujk{j}{k} = \left|\Im\left(\frac{\plgt(\alphajk{j}{k})}{\lambda_ka_k(\alphajk{j}{k})\plgt(\gamma)} \right)\right|
\end{equation}
which are independent of the choice of $\alphajk{j}{k}$. We are now ready to state the crucial technical lemma:

\begin{general}{Lemma}\label{lemmasignrelations}
The signs $\eps_j\djk{j}{k}$ obey the following relations:
\begin{itemize}
\item[$(\mathrm{i})$] $ \forall\, l \leq j, j' \leq n, \forall\, 1 \leq k < l, \eps_j \djk{j}{k} = \eps_{j'} \djk{j'}{k} $. In other words, the rows $1 \leq k < l$ in the sign table are constant and can be ignored. 
\item[$(\mathrm{ii})$] $\forall\, l \leq k, j \leq n, k \neq j, \eps_j\djk{j}{k} = -\sign(\lambda_j\lambda_k)\eps_k\djk{k}{j}$. In other words, the sign table is completely determined by its upper triangular portion.
\item[$(\mathrm{iii})$] $ \forall\, l \leq j, k, k' \leq n$, $k \neq k' \neq j$, if $\ujk{j}{k} \leq \ujk{j}{k'}$ then $\eps_{k}\djk{k}{k'} = \eps_{j}\djk{j}{k'}$. This technical point says that the knowledge of a single column is enough to determine the entire sign table, the particular sign relations being given by  the relative positions of the $\ujk{j}{k}$ for $k \neq j$ on the real axis.
\end{itemize}
\end{general}

\begin{proof}
(i) To prove the first relation, one only needs to notice that when $1 \leq k < l$ it is possible to choose $\alphajk{j}{k} = \alphajk{n}{k}$ for any $l \leq j \leq n$ which leads to $\eps_j\djk{j}{k} = \eps_n \djk{n}{k}$. The desired relation follows. \smallskip

(ii) To prove the second relation we remark that when $k \neq j$, $a_j(\alphajk{j}{k}) = -\lambda_k a_k(\alphajk{j}{k})/\lambda_j$ and therefore $a_j(-\sign(\lambda_j\lambda_k)\alphajk{j}{k}) > 0$. Since for all $k' \neq j, k$, $a_{k'}(\alphajk{j}{k}) = 0$, we may replace $\alphajk{k}{j}$ in the positive dual family to $a_1, \dots, \omitvar{a_k}, \dots, a_n$ by $-\sign(\lambda_j\lambda_k)\alphajk{j}{k}$. Thus, the sign $\djk{k}{j}$ satisfies by definition :
$$\djk{k}{j}\frac{\plgt(-\sign(\lambda_j\lambda_k)\alphajk{j}{k})}{\eps_k x(\gamma)} \in \hh $$
which we compare to
$$ \djk{j}{k} \frac{\plgt(\alphajk{j}{k})}{\eps_j \plgt(\gamma)} \in \hh$$
This gives exactly the relation $\eps_j\djk{j}{k} = -\sign(\lambda_j\lambda_k)\eps_k\djk{k}{j}$. \smallskip

(iii) This last point is more subtle and will only be used in section \ref{modextendpartfive}. For any three distinct indices $l \leq j, k, k' \leq n$, let us set:
$$\vjkl = \lambda_k^2 a_k(\alphajk{j}{k}) \alphajk{j}{k'} - \lambda_k\lambda_{k'} a_{k'}(\alphajk{j}{k'}) \alphajk{j}{k}$$
Then clearly $a_{i}(\vjkl) = 0$ whenever $i \neq j, k, k'$ and, by construction, it is also true that $a_{j}(\vjkl) = 0$. Indeed, using the fact that $\sum_{i = l}^n \lambda_i a_i = 0$ with $\lambda_{j} \neq 0$ we get:
\begin{align*}
\lambda_{j}a_{j}(\vjkl) & = \lambda_k^2\lambda_{j} a_k(\alphajk{j}{k}) a_{j}(\alphajk{j}{k'}) -\lambda_k \lambda_{k'}\lambda_{j} a_{k'}(\alphajk{j}{k'}) a_{j}(\alphajk{j}{k})\\
\lambda_{j}a_{j}(\vjkl) & = -\lambda_k^2 a_k(\alphajk{j}{k}) \sum_{i \neq j} \lambda_{i}a_{i}(\alphajk{j}{k'}) + \lambda_k\lambda_{k'} a_{k'}(\alphajk{j}{k'}) \sum_{i \neq j} \lambda_{i}a_{i}(\alphajk{j}{k})\\
\lambda_{j}a_{j}(\vjkl) & = -\lambda_k^2\lambda_{k'} a_k(\alphajk{j}{k}) a_{k'}(\alphajk{j}{k'}) + \lambda_k^2\lambda_{k'} a_{k'}(\alphajk{j}{k'}) a_{k}(\alphajk{j}{k})\\
\lambda_{j}a_{j}(\vjkl) & = 0.
\end{align*}
Moreover, $a_{k'}(\vjkl) = \lambda_{k}^2 a_{k}(\alphajk{j}{k}) a_{k'}(\alphajk{j}{k'}) > 0$. Thus, a possible choice for $\alphajk{k}{k'}$ is $\alphajk{k}{k'} = \vjkl$ and $x(\vjkl) \in \cc - \rr$. It follows that the sign $\djk{k}{k'}$ satisfies:
$$\djk{k}{k'}\frac{\plgt(\vjkl)}{\eps_{k} \plgt(\gamma)} \in \hh.$$
Replacing $\vjkl$ by its expression gives:
$$\djk{k}{k'}\left(\frac{\lambda_k^2 a_k(\alphajk{j}{k})\plgt(\alphajk{j}{k'})}{\eps_{k}\plgt(\gamma)} - \frac{\lambda_k\lambda_{k'} a_{k'}(\alphajk{j}{k'})\plgt(\alphajk{j}{k})}{\eps_{k}\plgt(\gamma)}\right) \in \hh $$
We may now rewrite this expression in terms of elements related to the values $\ujk{j}{k}$ and $\ujk{j}{k'}$:
\begin{equation}\label{equjk}
\eps_{k}\lambda_{k}^2\lambda_{k'}\djk{k}{k'}a_{k}(\alphajk{j}{k}) a_{k'}(\alphajk{j}{k'})\left(\frac{\plgt(\alphajk{j}{k'})}{\lambda_{k'}a_{k'}(\alphajk{j}{k'})\plgt(\gamma)} - \frac{\plgt(\alphajk{j}{k})}{\lambda_ka_k(\alphajk{j}{k})\plgt(\gamma)}\right) \in \hh.
\end{equation}
Let us denote by $\Ujk{j}{k}$ the complex number 
$$\Ujk{j}{k} = \frac{\plgt(\alphajk{j}{k})}{\lambda_ka_k(\alphajk{j}{k})\plgt(\gamma)}$$
so that \refp{equjk} may be rewritten as:
$$ \eps_{k}\lambda_{k'}\djk{k}{k'}(\Ujk{j}{k'} - \Ujk{j}{k}) \in \hh.$$
By definition $\ujk{j}{k} = |\Im(\Ujk{j}{k})|$ and $\ujk{j}{k'} = |\Im(\Ujk{j}{k'})|$,
thus the sign of $\eps_{k}\lambda_{k'}\djk{k}{k'}$ depends on which of the two values $\ujk{j}{k}, \ujk{j}{k'}$ is the largest. If $\ujk{j}{k} \leq \ujk{j}{k'}$ then the sign $\djk{k}{k'}$ satisfies:
$$\djk{k}{k'}\frac{\plgt(\alphajk{j}{k'})}{a_{k'}(\alphajk{j}{k'})\eps_{k}\plgt(\gamma)} \in \hh $$
which we compare to the sign $\djk{j}{k'}$ satisfying:
$$\djk{j}{k'}\frac{\plgt(\alphajk{j}{k'})}{a_{k'}(\alphajk{j}{k'})\eps_{j}\plgt(\gamma)} \in \hh.$$
This leads to the desired equality $\eps_k \djk{k}{k'} = \eps_j \djk{j}{k'}$. Note that if $\ujk{j}{k} = \ujk{j}{k'}$ then the two complex numbers $\Ujk{j}{k}$ and $\Ujk{j}{k'}$ must lie in opposite half-planes in $\cc-\rr$ and the sign equality holds.
\end{proof}

We are now ready to prove as a corollary that any intersection of three of the cones $\Cun{j}$ is empty. In fact, we prove something slightly stronger using the following definition.

\begin{general}{Definition}
Suppose that $l \leq j, j' \leq n$ are two distinct indices. We shall say that the two cones $\Cun{j}$ and $\Cun{j'}$ are compatible if the columns associated to the indices $j$ and $j'$ in the sign table are compatible. Explicitly, $\Cun{j}$ and $\Cun{j'}$ are compatible if and only if for all $1 \leq k \leq n$ (or $l \leq k \leq n$ by \lemmasignstoprows), $k \neq j, j'$ implies $ \eps_j\djk{j}{k} = \eps_{j'}\djk{j'}{k}$.
\end{general}

One may think of this definition as providing a necessary (but not sufficient) condition for two cones $\Cun{j}$ and $\Cun{j'}$ to have a non-empty intersection. We may now prove the lemma:

\begin{general}{Lemma}\label{lemmatroisvide}
\begin{itemize}
\item[$(\mathrm{i})$] For any two distinct indices $l \leq j, j' \leq n$, if $\Cun{j} \cap \Cun{j'} \neq \emptyset$ then the cones $\Cun{j}$ and $\Cun{j'}$ are compatible.
\item [$(\mathrm{ii})$] For any three distinct indices $l \leq j, j', j'' \leq n$, if $\eps_j\djk{j}{j''} = \eps_{j'}\djk{j'}{j''}$ and $\eps_{j'}\djk{j'}{j} = \eps_{j''}\djk{j''}{j}$ then $\eps_j\djk{j}{j'} = -\eps_{j''}\djk{j''}{j'}$.
\item [$(\mathrm{iii})$] For any three distinct indices $l \leq j, j', j'' \leq n$, if on the one hand the cones $\Cun{j}$ and $\Cun{j'}$ are compatible and on the other hand the cones $\Cun{j'}$ and $\Cun{j''}$ are compatible, then the cones $\Cun{j}$ and $\Cun{j''}$ are not compatible.
\item[$(\mathrm{iv})$] For any three distinct indices $l \leq j, j', j'' \leq n$, $\Cun{j} \cap \Cun{j'} \cap \Cun{j''} = \emptyset$.
\end{itemize}
\end{general}

\begin{proof}
(i) Suppose that the two distinct indices $l \leq j, j' \leq n$ are such that $\Cun{j}$ and $\Cun{j'}$ are not compatible. This means that there is an index $l \leq k \leq n$ distinct from $j, j'$ such that $\eps_j\djk{j}{k} = -\eps_{j'}\djk{j'}{k}$. Assume without loss of generality that $\eps_j\djk{j}{k} = +1$. Suppose that there is a vector $\delta \in \Cun{j} \cap \Cun{j'}$. Then on the one hand, $a_k(\delta) >0$ since $\eps_j\djk{j}{k} = +1$ and $\delta \in \Cun{j}$, while on the other hand, $a_k(\delta) \leq 0$ since $\eps_{j'}\djk{j'}{k} = -1$ and $\delta \in \Cun{j'}$. This yields a contradiction. \smallskip

(ii) Consider three indices $l \leq j \neq j' \neq j'' \leq n$. Using lemma \ref{lemmasignrelations}, (i) we get the following three relations:
\begin{align}
\eps_j\djk{j}{j'} &= -\sign(\lambda_j\lambda_{j'})\eps_{j'}\djk{j'}{j}\label{signone}\\
\eps_j\djk{j}{j''} &= -\sign(\lambda_j\lambda_{j''})\eps_{j''}\djk{j''}{j}\label{signtwo}\\
\eps_{j'}\djk{j'}{j''} &= -\sign(\lambda_{j'}\lambda_{j''})\eps_{j''}\djk{j''}{j'}. \label{signthree}
\end{align}
Suppose that:
\begin{align}
\eps_j\djk{j}{j''} &= \eps_{j'}\djk{j'}{j''}\label{signfour}\\
\eps_{j'}\djk{j'}{j}& = \eps_{j''}\djk{j''}{j}\label{signfive}
\end{align}
We wish to prove that $\eps_{j}\djk{j}{j'} = -\eps_{j''}\djk{j''}{j'}$. We start by rewriting \refp{signone} using \refp{signfive} as:
\begin{align*} 
\eps_j\djk{j}{j'} &= -\sign(\lambda_j\lambda_{j'})\eps_{j''}\djk{j''}{j}\\
\eps_j\djk{j}{j'} &= \sign(\lambda_j\lambda_{j'})\sign(\lambda_j\lambda_{j''})\eps_j\djk{j}{j''} \\
\eps_j\djk{j}{j'} &= \sign(\lambda_j^2\lambda_{j'}\lambda_{j''})\eps_{j'}\djk{j'}{j''} \\
\eps_j\djk{j}{j'} &= -\sign(\lambda_j^2\lambda_{j'}^2\lambda_{j''}^2)\eps_{j''}\djk{j''}{j'}
\end{align*}
where we used \refp{signtwo} in the second line, \refp{signfour} in the first line and \refp{signfive} in the final line. Thus $\eps_j\djk{j}{j'} = -\eps_{j''}\djk{j''}{j'}$. \smallskip

(iii) Consider three indices $l \leq j \neq j' \neq j'' \leq n$ and suppose that $\Cun{j}$ is compatible with $\Cun{j'}$ while $\Cun{j'}$ is compatible with $\Cun{j''}$. Then in particular $\eps_j\djk{j}{j''} = \eps_{j'}\djk{j'}{j''}$ and $\eps_{j'}\djk{j'}{j} = \eps_{j''}\djk{j''}{j}$ therefore by (ii) we get $\eps_j\djk{j}{j'} = -\eps_{j''}\djk{j''}{j'}$ which guarantees that $\Cun{j}$ and $\Cun{j''}$ are incompatible. \smallskip

(iv) Consider three indices $l \leq j \neq j' \neq j'' \leq n$ and suppose that $\Cun{j} \cap \Cun{j'} \cap \Cun{j''} \neq \emptyset$. In particular, each of the two-ways intersections are non-empty, and therefore by (i), the of cones $\Cun{j}$ and $\Cun{j'})$ are compatible, and the same is true for the cones $\Cun{j'}$ and $\Cun{j''}$, as well as for the cones $\Cun{j''}$ and $\Cun{j})$. This contradicts (iii), therefore we must have an empty intersection $\Cun{j} \cap \Cun{j'} \cap \Cun{j''} = \emptyset$.
\end{proof}

It follows from lemma \ref{lemmatroisvide}, (iv) that any vector $\delta \in V$ belongs to at most two of the cones $\Cun{j}$. The fourth and next step in the proof of Proposition \ref{propmodularextend} is showing that a vector $\delta \in V$ cannot belong to exactly one of the cones $\Cun{j}$.

\subsubsection{Any vector $\delta \in V$ is in none or exactly two of the cones $\Cun{j}$}\label{modextendpartfive}

In this section, we prove that a vector $\delta \in V$ belongs to either either $0$ or $2$ of the $\Cun{j}$'s. From lemma \ref{lemmatroisvide}, (iv) we already know that a vector $\delta$ belongs to at most $2$ of the cones $\Cun{j}$, thus it will suffice to prove that $\delta$ cannot belong to exactly one cone $\Cun{j}$. The result will be proven in the following form:

\begin{general}{Lemma}\label{lemmastepfive}
For all $l \leq j \leq n$:
$$\Cun{j} \subset \bigcup_{k \neq j} \Cun{k}$$
where the union ranges on indices $l \leq k \leq n$, $k \neq j$. 
\end{general}

This is the most technical part of the proof and it will use most results of section \ref{modextendpartthree}. It will follow from this result that if $\delta \in V$ belongs to the cone $\Cun{j}$ then, since $\Cun{j} \subset \cup_{k \neq j} \Cun{k}$, there is an index $k \neq j$ such that $\delta \in \Cun{k}$. Hence $\delta$ cannot belong to exactly one of the cones $\Cun{j}$.
Let us now sketch the proof of Lemma \ref{lemmastepfive}. Define for all $l \leq j \leq n$ the set of indices: 
\begin{equation}\label{ijfinerlemmastepfive}
\mathcal{I}(j) = \{ l \leq j' \leq n \setseparator j' \neq j \text{ and the cones } \Cun{j}, \Cun{j'} \text{ are compatible}\}.
\end{equation}
From Lemma \ref{lemmatroisvide} it is clear that we only need to prove that for all $l \leq j \leq n$:
\begin{equation}\label{finerlemmastepfive}
\Cun{j} \subset \bigcup_{j' \in \mathcal{I}(j)} \Cun{j'}.
\end{equation}
To do so, we shall first prove \refp{finerlemmastepfive} when $\cardinalshort{\mathcal{I}(j)} \geq 2$. Then, we prove that for any $l \leq j \leq n$, $\cardinalshort{\mathcal{I}(j)} \geq 1$, and lastly, we prove \refp{finerlemmastepfive} when $\cardinalshort{\mathcal{I}(j)} = 1$. The hypothesis that $k^{-}(a_1, \dots, a_n) > 0$ will only be used in this last step. 

Let us now prove that for any $l \leq j \leq n$, if $\cardinalshort{\mathcal{I}(j)} \geq 2$ then \refp{finerlemmastepfive} holds. To this end we introduce the following partition of the cones $\Cun{j}$. Write as before $\hjplus{j} = \{v \in V \setseparator a_j(v) > 0
\}$ and $\hjminus{j} = V - \hjplus{j}$. Then the two cones $\Cun{j} \cap \hjplus{j}$
and $\Cun{j} \cap \hjminus{j}$ form a partition of $\Cun{j}$. Our claim is that whenever some $\Cun{j}$ intersects some other $\Cun{j'}$ for $j \neq j'$, the intersection must be exactly one of these components. More precisely:

\begin{general}{Lemma}\label{lemmacomponent}
\begin{itemize}
\item[$(\mathrm{i})$] Suppose that the cones $\Cun{j}$ and $\Cun{j'}$ are compatible for some distinct indices $l \leq j, j' \leq n$. Then:
$$\Cun{j} \cap \Cun{j'} = \Cun{j} \cap \hj{\eps_{j'}\djk{j'}{j}}{j} = \Cun{j'} \cap \hj{\eps_j\djk{j}{j'}}{j'}.$$
\item[$(\mathrm{ii})$] As a consequence, for any $l \leq j \leq n$, if $\mathcal{I}(j)$ contains at least two distinct indices $l \leq j', j'' \leq n$, $j \neq j', j''$ then $\Cun{j} \subset \Cun{j'} \cup \Cun{j''}$.
\end{itemize}
\end{general} 

\begin{proof}
(i) The cones $\Cun{j}$ and $\Cun{j'}$ are compatible, so by definition we get that for all index $1 \leq k\leq n$ distinct from $j$ and $j'$, $\eps_{j}\djk{j}{k} = \eps_{j'}\djk{j'}{k}$. Recall that the half-spaces $H_k^{\pm}$ are defined for $1 \leq k \leq n$ by $H_k^{+} = \{\delta \in V \setseparator a_k(\delta) > 0\}$ and $H_k^{-} = V - H_k^{+}$. For simplicity, let us denote by $H_k$ the half-space $H_k^{+}$ if $\eps_j\djk{j}{k} = 1$ or $H_k^{-}$ if $\eps_j \djk{j}{k} = -1$ for $1 \leq k \leq n$, $k \neq j$. Denote also by $H_j$ the half-space $H_j^{+}$ if $\eps_{j'}\djk{j'}{j} = 1$ or $H_j^{-}$ if $\eps_{j'} \djk{j'}{j} = -1$. Then
$$ \Cun{j} = \bigcap_{k \neq j} H_k~~~~\Cun{j} = \bigcap_{k \neq j'} H_k$$
so that 
$$\Cun{j} \cap \Cun{j'} = \left(\bigcap_{k \neq j, j'} H_k\right) \cap H_j \cap H_{j'} =  \Cun{j} \cap H_j = \Cun{j'} \cap H_{j'}$$
which is the desired result.

(ii) Let us now suppose that $\mathcal{I}(j)$ contains at least two distinct indices $j', j''$. Then by definition the cones $\Cun{j}$ and $\Cun{j'}$ are compatible and so are the cones $\Cun{j}$ and $\Cun{j'}$. It follows from the Lemma \ref{lemmatroisvide}, (iii) that the cones $\Cun{j'}$ and $\Cun{j''}$ are not compatible, and more precisely, it follows from of Lemma \ref{lemmatroisvide}, (ii) that:
$$\eps_{j'}\djk{j'}{j} = - \eps_{j''}\djk{j''}{j} $$
Up to switching $j'$ and $j''$ we may assume without loss of generality that $\eps_{j'}\djk{j'}{j} = +1$ so that $\eps_{j''}\djk{j''}{j} = -1$. Then, using (i) we get:
\begin{align*}
\Cun{j} \cap \Cun{j'} & =\Cun{j} \cap \hjplus{j}\\
\Cun{j} \cap \Cun{j''} &= \Cun{j} \cap \hjminus{j}.
\end{align*}
Therefore:
$$\Cun{j} \cap (\Cun{j'} \cup \Cun{j''}) = \left(\Cun{j} \cap \Cun{j'}\right) \cup \left(\Cun{j} \cap \Cun{j''}\right) = (\Cun{j} \cap \hjplus{j}) \cup (\Cun{j} \cap \hjminus{j}) = \Cun{j} $$
and $\Cun{j} \subset \Cun{j'} \cup \Cun{j''}$. This proves \refp{finerlemmastepfive} in the case where $\cardinalshort{\mathcal{I}(j)} \geq 2$.
\end{proof}

The remainder of this section is devoted to the proof of \refp{finerlemmastepfive} in the case where $\mathcal{I}(j)$ contains no more than one element. In fact, we shall first prove using \lemmasignsujk that the set $\mathcal{I}(j)$ cannot be empty and then handle the case where $\mathcal{I}(j)$ contains exactly one element.

\begin{general}{Lemma}\label{lemmanonvide}
Fix an index $l \leq j \leq n$. Let $l \leq j' \leq n$ be an index satisfying $j' \neq j$ and for all $l \leq k \leq n$, $k \neq j, j'$, $\ujk{j}{j'} \leq \ujk{j}{k}$ (see \refp{defujk} for the definition of the positive real numbers $\ujk{j}{k}$).
Then the cones $\Cun{j}$ and $\Cun{j'}$ are compatible and $j' \in \mathcal{I}(j)$. In particular, $\mathcal{I}(j)$ cannot be empty.
\end{general}

\begin{proof}
By definition of $j'$ we get that for all indices $l \leq k \leq n$ distinct from both $j$ and $j'$, $\ujk{j}{j'} \leq \ujk{j}{k}$. Using \lemmasignsujk we get that for all such indices $l \leq k \leq n$ distinct from both $j$ and $j'$, $\eps_j\djk{j}{k} = \eps_{j'}\djk{j'}{k}$. Thus, $\Cun{j}$ and $\Cun{j'}$ are compatible and $j' \in \mathcal{I}(j)$.
\end{proof}

In the last part of the proof of lemma \ref{lemmastepfive} we need to treat the case where $\mathcal{I}(j)$ contains only one index $j'$. In that case, we wish to prove that one of the two components $\Cun{j} \cap \hjplus{j}$, $\Cun{j} \cap \hjminus{j}$ is empty. To achieve this, we first give a sufficient condition for one of these components to be empty under the condition that $k^{-}(a_1, \dots, a_n) > 0$.

\begin{general}{Lemma}\label{lemmasuffisantvide}
Fix a sign $\nu \in \{-1, 1\}$ and an index $l \leq j \leq n$. Assume that for all index $l \leq k \leq n$ distinct from $j$, $\nu\lambda_j \lambda_k \eps_j \djk{j}{k} > 0$. Then $\Cun{j} \cap \hj{\nu}{j} = \emptyset$.
\end{general}

\begin{proof}
Here we use the hypothesis that $k^{-}(a_1, \dots, a_n) > 0$. Recall that $\sum_{k = l}^n \lambda_k a_k = 0$ is the standard non-trivial relation among $a_1, \dots, a_n$ with the coefficients $\lambda_k$ satisfying \refp{signslambda}. In particular, $\lambda_l < 0$ and $\lambda_n > 0$. Suppose that $\Cun{j} \cap \hj{\nu}{j} \neq \emptyset$ and fix some $\delta \in \Cun{j} \cap \hj{\nu}{j}$. By definition of $\Cun{j}$, for all index $1 \leq k \leq n$ distinct from $j$, the vector $\delta$ satisfies $a_k(\delta)  > 0$ if $\eps_j \djk{j}{k} = 1$ and $a_k(\delta) \leq 0$ if $\eps_j \djk{j}{k} = -1$. In addition, by definition of $\hj{\nu}{j}$, 
this element $\delta$ must satisfy $a_j(\delta) > 0$ if $\nu = 1$ and $a_j(\delta) \leq 0$ if $\nu = -1$. In particular, $\nu \lambda_j^2 a_j(\delta) \geq 0$ in either case. Thus, the assumption that $\nu \lambda_j \lambda_k \eps_j\djk{j}{k} > 0$ for all $l \leq k \leq n$, $k \neq j$ leads to $\nu \lambda_j \lambda_k a_k(\delta) \geq 0$ for all $l \leq k \leq n$. Using the relation $\sum_{k =l}^n \lambda_k a_k = 0$ we get the equality:
$$\sum_{k = l}^n \nu \lambda_j \lambda_k a_k(\delta) =0 $$
which is the vanishing of a sum of non-negative terms. Therefore, for each $l \leq k \leq n$ we get $\nu \lambda_j \lambda_k a_k(\delta) = 0$ and thus $a_k(\delta) = 0$. From the vanishing of $a_j(\delta)$ we obtain that $\nu = -1$ while the vanishing of $a_k(\delta)$ for $k \neq j$ leads to $\eps_j \djk{j}{k} = -1$ for $l \leq k \leq n$, $k \neq j$. Thus the assumption $\nu \lambda_j \lambda_k \eps_j\djk{j}{k} > 0$ for $l \leq k \leq n$, $k \neq j $ gives $\lambda_j \lambda_k > 0$ for all $l \leq k \leq n$ and all the coefficients $\lambda_l, \dots, \lambda_n$ must share the same sign, which contradicts the fact that $\lambda_l < 0$ and $\lambda_n > 0$. Therefore we must conclude that $\Cun{j} \cap \hj{\nu}{j} = \emptyset$ as claimed.
\end{proof}

Using lemma \ref{lemmasuffisantvide} we shall finally prove that if $\mathcal{I}(j)$ contains only one index $j'$, then $\Cun{j} \subset \Cun{j'}$. This will complete the proof of lemma \ref{lemmastepfive}.

\begin{general}{Lemma}\label{finallemma}
Suppose that the index $l \leq j \leq n$ is such that $\mathcal{I}(j)$ contains exactly one index $j'$. Then:
\begin{itemize}
\item[$(\mathrm{i})$] For all $l \leq k \leq n$, $k \neq j$, the sign $\sign(\lambda_k)\djk{j}{k}$ is equal to $\sign(\lambda_{j'})\djk{j}{j'}$.
\item[$(\mathrm{ii})$] As a consequence, $\Cun{j} \subset \Cun{j'}$. 
\end{itemize}
\end{general}

\begin{proof}
(i) Let us fix a bijection $\sigma : \{1, \dots, n-l\} \to \{l \leq k \leq n \setseparator k \neq j\}$ such that $\ujk{j}{\sigma(k)} \leq \ujk{j}{\sigma(k')}$ whenever $1 \leq k < k' \leq n-l$. In particular, it follows from lemma \ref{lemmanonvide} that $j' = \sigma(1)$. Let us then prove by induction on $2 \leq \kappa \leq n-l$ that $\sign(\lambda_{\sigma(\kappa)}) \djk{j}{\sigma(\kappa)} = \sign(\lambda_{j'})\djk{j}{j'}$. 

\underline{First case:} $\kappa = 2$. We set $k = \sigma(\kappa)$. By definition of $\sigma$ we get that for $3 \leq \kappa' \leq n-l$, $\ujk{j}{k} \leq \ujk{j}{\sigma(\kappa')}$. From \lemmasignsujk we get $\eps_j \djk{j}{\sigma(\kappa')} = \eps_k \djk{k}{\sigma_j(\kappa')}$ for any $3 \leq \kappa' \leq n-l$. Thus $k \not\in \mathcal{I}(j)$ if and only if $\eps_j\djk{j}{j'} = - \eps_k \djk{k}{j'}$ where once again $j' = \sigma(1)$. Using \lemmasignstranspose this gives:
$$\eps_j\djk{j}{j'} = \sign(\lambda_{j'}\lambda_k) \eps_{j'} \djk{j'}{k} = \sign(\lambda_{j'}\lambda_k) \eps_{j} \djk{j}{k}.$$
where the last equality holds as the cones $\Cun{j}$ and $\Cun{j'}$ are compatible by assumption. Hence, $\sign(\lambda_{j'})\djk{j}{j'} = \sign(\lambda_k) \djk{j}{k}$ and the case $\kappa = 2$ is proven.

\underline{Induction:} Assume that the result holds for $k = \sigma(2), \dots, \sigma(\kappa-1)$. Let $k = \sigma(\kappa)$. For any $\kappa < \kappa' \leq n-l$, by definition of $\sigma$ and \lemmasignsujk we get $\eps_j \djk{j}{\sigma_j(k')} = \eps_k \djk{k}{\sigma_j(k')}$. Therefore, $k \not\in \mathcal{I}_j$ is equivalent to the existence of some $k'' = \sigma(\kappa'')$ with $\kappa'' < \kappa$ such that $\eps_j\djk{j}{k''} = - \eps_k \djk{k}{k''}$. Using once again \lemmasignstranspose we rewrite this as: 
$$\eps_j\djk{j}{k''} = \sign(\lambda_k \lambda_{k''}) \eps_{k''} \djk{k''}{k} =  \sign(\lambda_k \lambda_{k''}) \eps_{j} \djk{j}{k}$$
where the last equality holds by \lemmasignsujk as $\ujk{j}{k''} \leq \ujk{j}{k}$. Thus we get $\sign(\lambda_{k''})\djk{j}{k''} = \sign(\lambda_k) \djk{j}{k}$. It then follows from the induction hypothesis $\sign(\lambda_{k''})\djk{j}{k''} = \sign(\lambda_{j'})\djk{j}{j'}$ that $\sign(\lambda_{j'})\djk{j}{j'} = \sign(\lambda_k)\djk{j}{k}$. This completes the proof by induction.

(ii) Let $\nu = \sign(\lambda_j \lambda_{j'}) \eps_j \djk{j}{j'}$. For all index $l \leq k\leq n$ distinct from $j$, since by (i) $\sign(\lambda_k)\djk{j}{k} = \sign(\lambda_{j'}) \djk{j}{j'}$ we get:
\begin{align*}
\nu \lambda_j \lambda_k \eps_j \djk{j}{k} & = |\lambda_j| \sign(\lambda_{j'})\djk{j}{j'}\lambda_k \eps_j^2\djk{j}{k} \\
\nu \lambda_j \lambda_k \eps_j \djk{j}{k} & = |\lambda_j|\sign(\lambda_k)\lambda_k (\djk{j}{k})^2 \\
\nu \lambda_j \lambda_k \eps_j \djk{j}{k} & = |\lambda_j\lambda_k| > 0.
\end{align*}
It then follows from lemma \ref{lemmasuffisantvide} that $\Cun{j} \cap \hj{\nu}{j} = \emptyset$ and since $\Cun{j} = (\Cun{j} \cap \hj{-\nu}{j}) \cup (\Cun{j} \cap \hj{\nu}{j})$ we conclude that $\Cun{j} = \Cun{j} \cap \hj{-\nu}{j}$. From lemma \ref{lemmasignrelations}, (ii) we get $\nu = -\eps_{j'}\djk{j'}{j}$, therefore $\Cun{j} = \Cun{j} \cap \hj{-\nu}{j} = \Cun{j} \cap \Cun{j'}$ by lemma \ref{lemmacomponent}. This gives the conclusion $\Cun{j} \subset \Cun{j'}$.
\end{proof}

Let us now piece all these results together to prove lemma \ref{lemmastepfive}. \bigskip

\begin{proofbis}{Proof of lemma \ref{lemmastepfive}}
Consider $l \leq j \leq n$. Then $\cardinalshort{\mathcal{I}}(j) \geq 1$ by lemma \ref{lemmanonvide}. If $\cardinalshort{\mathcal{I}(j)} = 1$ then $\mathcal{I}(j)$ contains exactly one element $j' \neq j$ and $\Cun{j} \subset \Cun{j'}$ by lemma \ref{finallemma}, (ii). If $\cardinalshort{\mathcal{I}(j)} \geq 2$ then $\mathcal{I}(j)$ contains at least two distinct elements $j', j'' \neq j$ and $\Cun{j} \subset \Cun{j'} \cup \Cun{j''}$ by lemma \ref{lemmacomponent}, (ii). Thus in any case, $\Cun{j} \subset \cup_{k \neq j} \Cun{k}$ and the proof is complete.
\end{proofbis}

\subsubsection{Vanishing of $\fun$ and proof of Proposition \ref{propmodularextend}}\label{modextendpartfour}

In this section we prove that $\fun$ vanishes. From the previous sections, we already know that any vector $\delta \in V$ belongs to either $0$ or $2$ of the cones $\Cun{j}$. In the case where $\delta$ belongs to none of the $\Cun{j}$'s, it is clear that $\fun(\delta) = 0$ by definition. Thus, it will suffice to prove that when $\delta$ belongs to the intersection of any two cones $\Cun{j}$ and $\Cun{j'}$ for $j \neq j'$, $\fun(\delta) = 0$.

\begin{general}{Lemma}\label{lemmafunzero}
Assume $\delta \in C^{1}_j \cap C^{1}_{j'}$ for some distinct indices $l \leq j, j' \leq n$. Then $\fun(\delta) = 0$.
\end{general}

\begin{proof}
It follows from lemma \ref{lemmatroisvide}, (iv) that for any $j'' \neq j, j'$, $\delta \not\in C^{1}_{j''}$. Thus by definition of $\fun$ (see lemma \ref{lemmafunfdeux}):
$$ \fun(\delta) = (-1)^{j + 1 + D_j + n\frac{1+\eps_j}{2}} + (-1)^{j'+1 + D_{j'} + n\frac{1+\eps_{j'}}{2}}.$$ 
Let $g(j, j') = (j-j') + (D_j - D_{j'}) + n\frac{\eps_j - \eps_{j'}}{2}$. The statement that $\fun(\delta) = 0$ is equivalent to $g(j, j') \equiv 1 \mod 2$ as 
$$\fun(\delta) = (-1)^{j+1+D_j + n\frac{1-\eps_j}{2}}(1 + (-1)^{g(j, j')}).$$
Let us then prove that $g(j, j') \equiv 1 \mod 2$. Write $D_j$ explicitly as:
$$D_j = \sum_{k \neq j} \frac{\djk{j}{k} -1}{2}.$$
Since $\Cun{j} \cap \Cun{j'} \neq \emptyset$, it follows from lemma \ref{lemmatroisvide}, (i) that for all index $1 \leq k \leq n$ distinct from both $j$ and $j'$, $\eps_j \djk{j}{k} = \eps_{j'} \djk{j'}{k}$. Thus:
\begin{align*}
D_j - D_{j'} &= \sum_{k \neq j} \frac{\djk{j}{k} -1}{2} - \sum_{k \neq j'} \frac{\djk{j'}{k} -1}{2} \\ 
D_j - D_{j'} &= \frac{\djk{j}{j'} - \djk{j'}{j}}{2}+ \sum_{k \neq j, j'} \frac{\djk{j}{k} -\djk{j'}{k}}{2}   \\
D_j - D_{j'} &= \frac{\djk{j}{j'} - \djk{j'}{j}}{2} + \frac{(1-\eps_j\eps_{j'})}{2}\sum_{k \neq j, j'} \djk{j}{k}.
\end{align*}
The sum over $k \neq j, j'$ contains $n-2$ terms which are either $+1$ or $-1$ therefore we already obtain
$$D_j - D_{j'} \equiv  \frac{\djk{j}{j'} - \djk{j'}{j}}{2} + (n-2) \frac{1-\eps_j \eps_{j'}}{2} \equiv  \frac{\djk{j}{j'} - \djk{j'}{j}}{2} + (n-2) \frac{\eps_j-\eps_{j'}}{2} \mod 2$$
and therefore
\begin{align*}
g(j, j') &\equiv (j - j') + \frac{\djk{j}{j'} - \djk{j'}{j}}{2} + (2n-2) \frac{\eps_j -\eps_{j'}}{2} \mod 2\\
g(j, j') &\equiv (j - j') +\frac{\djk{j}{j'} - \djk{j'}{j}}{2}\mod 2.
\end{align*}
Recall that lemma \ref{lemmasignrelations}, (ii) gives the equality $\eps_j \djk{j}{j'} = -\sign(\lambda_j\lambda_{j'}) \eps_{j'}\djk{j'}{j}$ and that by definition $\eps_j = (-1)^{j+n} \sign(\lambda_j)$. Thus, $(-1)^{j}\djk{j}{j'} = -(-1)^{j'} \djk{j'}{j}$ and:
$$g(j, j') \equiv (j-j') + \djk{j}{j'}\frac{1+(-1)^{j+j'}}{2} \mod 2.$$
There are only two cases to treat, depending on the parity of $j-j'$. First, suppose that $j-j'$ is even. Then $1+ (-1)^{j+j'} = 2$ and $g(j, j') \equiv \djk{j}{j'} \equiv 1 \mod 2$. Now, suppose that $j-j'$ is odd. Then $1+(-1)^{j+j'} = 0$ and $g(j, j') \equiv j - j' \equiv 1 \mod 2$. In each of these two cases, $g(j, j') \equiv 1 \mod 2$ and thus $\fun(\delta) = 0$. This completes the proof.
\end{proof}

Finally, we may piece together all the results from sections \ref{modextendpartone} to \ref{modextendpartfour} to give the proof of Proposition \ref{propmodularextend} and thus deduce Theorem \ref{theoremmodularextend}. \medskip

\begin{proofbis}{Proof of Proposition \ref{propmodularextend}}
Let $a_1, \dots, a_n \in \Lambda$ be primitive integral linear forms on $V$ satisfying $\rank(a_1, \dots, a_n) = n-1$ and such that $0$ is not a barycenter of $a_1, \dots, a_n$. This means that there is a relation $\sum_{j = 1}^n \lambda_j a_j = 0$ with at least one positive and one negative coefficient among the $\lambda_j$'s. Using the results from sections \ref{modextendpartone} and \ref{modextendparttwo} we may suppose that the coefficients $\lambda_1, \dots, \lambda_n$ in the standard non-trivial relation among $a_1, \dots, a_n$ satisfy the relations \refp{signslambda} and that
$$\prod_{j = l}^{n} G_{n-2, a_1, \dots, \omitvar{a_j}, \dots, a_n}(v)(w, \plgt)^{(-1)^{j+1}} = \prod_{\delta \in (v+L)/\zz\gamma}\left(1-e^{2i\pi\left(\frac{w+\plgt(\delta)}{\plgt(\gamma)}\right)}\right)^{\fun(\delta)} \left(1-e^{-2i\pi\left(\frac{w+\plgt(\delta)}{\plgt(\gamma)}\right)}\right)^{\fdeux(\delta)}.$$
Now, suppose that $\fun \neq 0$. Then there exists some $\delta \in V$ such that $\fun(\delta) \neq 0$. Therefore, by definition of $\fun$, there is some index $l \leq j \leq n$ such that $\delta \in \Cun{j}$. By lemma \ref{lemmastepfive} the cone $\Cun{j}$ is a subset of $\cup_{k \neq j} \Cun{k}$ and therefore there is an index $l \leq j' \leq n$ distinct from $j$ such that $\delta \in \Cun{j'}$. Then, since $\delta \in \Cun{j} \cap \Cun{j'}$, using lemma \ref{lemmafunzero} we get that $\fun(\delta) = 0$ which is a contradiction. Therefore $\fun = 0$ and similarly $\fdeux = 0$, which gives the desired conclusion:
$$\prod_{j = 1}^n G_{n-2, a_1, \dots, \omitvar{a_j}, \dots, a_n}(v)(w,x)^{(-1)^{j+1}} = 1.$$
\end{proofbis}

\textbf{Remark:} Note that to prove the similar case $\fdeux = 0$ one may reuse most of the work carried out in sections \ref{modextendpartone} to \ref{modextendpartfour}. Indeed, it is clear that the sign table governing the cones $\Cun{j}$ also governs the cones $\Cdeux{j}$ and that two cones $\Cdeux{j}$ and $\Cdeux{j'}$ are compatible if and only if $\Cun{j}$ and $\Cun{j'}$ are compatible. Thus, it is clear that lemma \ref{lemmatroisvide}, (iv) may be adapted for the cones $\Cdeux{j}$ as:
$$\textit{ For any three distinct indices } l \leq j, j', j'' \leq n, \Cdeux{j} \cap \Cdeux{j'} \cap \Cdeux{j''} = \emptyset.$$
Next, there is only a small adaptation to make to express a version of lemma \ref{lemmacomponent}, (i) for the cones $\Cdeux{j}$. Indeed, if $\Cdeux{j}$ and $\Cdeux{j'}$ are compatible then:
$$ \Cdeux{j} \cap \Cdeux{j'} = \Cdeux{j} \cap \hj{-\eps_{j'}\djk{j'}{j}}{j} = \Cdeux{j'} \cap \hj{-\eps_j\djk{j}{j'}}{j'}.$$
Another adaptation is in order for lemma \ref{lemmasuffisantvide} as the condition $\nu\lambda_j \lambda_k \eps_j \djk{j}{k} > 0$ for all index $l \leq k \leq n$ distinct from $j$ implies that $\Cun{j} \cap \hj{-\nu}{j} = \emptyset$.
Lemma \ref{lemmastepfive} may then be directly adapted as for all $l \leq j \leq n$:
$$\Cdeux{j} \subset \bigcup_{k \neq j}\Cdeux{k}.$$
Finally, a computation similar to the one carried out in the proof of lemma \ref{lemmafunzero} shows that $\fdeux = 0$.

\section{A classic smoothing operation}\label{sectionsmoothing}

In this section we are interested in the smoothed versions of both $G_{n-2, a_1, \dots, a_{n-1}}$ and $B_{n, a_1, \dots, a_n}$ functions which are defined by:
$$\smoothedgr = \frac{G_{n-2, a_1, \dots, a_{n-1}}(v)(w, x, L')^N}{G_{n-2, a_1, \dots, a_{n-1}}(v)(w, x, L)}.$$
and
$$\smoothedbn = NB_{n, a_1, \dots, a_n}(v)(w,x, L') - B_{n, a_1, \dots, a_n}(v)(w,x,L)$$
for linear forms $a_1, \dots, a_n$ which are primitive on both $L$ and $L'$ where $L/L' \simeq \zsz{N}$. In this section we prove Theorem \ref{theoremsmoothing} and explain how we deduce Theorem \ref{maintheorem}. Let us now give an overview of this section.  In section \ref{sectionsmoothinggeometricsetup} we give the geometric setup needed for the rest of the proof and give an explicit formula for the smoothed $\smoothedbn$ when the linear forms $a_1, \dots, a_n$ are linearly independent and the smoothing lattice $L'$ is \good for $a_1, \dots, a_n$. Then, in section \ref{sectionrational} we prove that the function $\smoothedbn$ is in fact a rational valued function which depends only on the linear forms $a_1, \dots, a_n$ and on the class of $v$ in $V/L'$ but not on $w, x \in \cc \times \homlc$. Finally, in section \ref{sectiontraces} we use Fourier analysis following \cite{CD} to prove that the rational numbers $\smoothedbn$ may be expressed in terms of traces of cyclotomic units from which we may deduce a bound on its denominator in terms of the dimension $n$ and the smoothing index $N$. At the end of this section, we shall explain how to derive Theorem \ref{maintheorem} from this last result.

\subsection{Geometric setup}\label{sectionsmoothinggeometricsetup}

In this section we consider the situation where $V$ is a $\qq$-vector space of dimension $n$ and $L$ is a lattice of rank $n$ in $V$ with a $\zz$-basis $B = [e_1, \dots, e_n]$. We fix an integer $N \geq 2$ and the smoothing lattice $L' = N\zz e_1 \oplus \zz e_2 \oplus \dots \oplus \zz e_n$. Denote by $\Lambda = \homlz$ and $\Lambda' = \homlprimez$ the dual spaces attached to $L$ and $L'$ respectively. Define $C = [f_1, \dots, f_n]$ the $\zz$-basis of $\Lambda$ dual to $B$ such that $\forall\, 1 \leq j,k \leq n, f_j(e_k) = \delta_{jk}$. Similarly, the $\zz$-basis $C'$ of $\Lambda'$ dual to $B'$ is given by $C' = [f_1/N, f_2, \dots, f_n]$. One may view $\Lambda \subset \Lambda'$ as rank $n$ lattices in the dual space $\dualsimple{V} = \mathrm{Hom}_{\qq}(V, \qq)$. Let us now define $\Lambda_N \subset \dualsimple{V}$ to be the set of linear forms $a \in \dualsimple{V}$ which restrict to primitive integral linear forms on both $L$ and $L'$. Explicitly:
$$ \Lambda_N = \left\{ \sum_{k = 1}^{n} \mu_kf_k \setseparator \mu_1, \dots, \mu_n \in \zz,\, \pregcd(\mu_1N, \mu_2, \dots, \mu_n) = 1\right\}.$$
The set $\Lambda_N$ is endowed with an action of the congruence subgroup $\congruencegrouphnn \subset \slnz{n}$ (see \refp{definitionmodulargroup}) given by $g\cdot(\mu_1, \mu_2, \dots, \mu_n) = (\mu_1, \mu_2, \dots, \mu_n) \times g^{-1}$. For the rest of this section we fix non-zero linear forms $a_1, \dots, a_n$ in $\Lambda_N$ which are linearly independent. This fixes a unique family of primitive vectors $\alpha_1, \dots, \alpha_n \in L$ such that for all $1 \leq j \leq n$:
$$a_j(\alpha_k) = 0, \forall\, 1 \leq k \neq j \leq n ~~ \text{ and } ~~ a_j(\alpha_j) = s_j > 0 $$
This is the primitive positive dual basis to $a_1, \dots, a_n$ in $L$ in the sense of [\hspace{1sp}\cite{firstpaper}, lemma 6]. This definition comes from the theory of rational polyhedral cones as the cone:
$$C := \{\delta \in V \setseparator \forall 1 \leq j \leq n, a_j(\delta) \geq 0 \}$$
can be expressed in terms of generators as:
$$C = \qq_{\geq 0} \alpha_1 + \dots + \qq_{\geq 0} \alpha_n.$$
We now recall the definition of a \good smoothing lattice.

\begin{general}{Definition}\label{defgoodsmoothing}
Suppose that $a_1, \dots, a_n \in \Lambda$ are linearly independent primitive linear forms and fix $\alpha_1, \dots, \alpha_n$ the primitive positive dual basis to $a_1, \dots, a_n$ in $L$. The smoothing lattice $L'$ of index $N$ in $L$ is said to be \good for $a_1, \dots, a_n$ if $a_1, \dots, a_n \in \Lambda_N$ and the primitive positive dual basis $\alpha'_1, \dots, \alpha'_n$ to $a_1, \dots, a_n$ in $L'$ is precisely $N\alpha_1, \dots, N\alpha_n$.
\end{general}

Notice that in general if $k_j$ is the order of $\alpha_j \mod L'$ in the cyclic group $L/L'$ then it is clear that $\alpha'_j = k_j\alpha_j$, thus the condition $\alpha_j \mod L'$ generates $L/L'$ is equivalent to $\alpha'_j = N\alpha_j$. In the case where $a_1, \dots, a_n$ are not linearly independent, we say that the smoothing lattice $L'$ is \good for $a_1, \dots, a_n$ if $a_1, \dots, a_n \in \Lambda_N$. Definition \refp{defgoodsmoothing} is inspired by the definition of a good smoothing ideal in \cite{DasguptaShintani}. Notice that in our case we do not suppose that $N$ is prime. Let us now give an equivalent formulation of this statement as a condition on the coordinates of $\alpha_1, \dots, \alpha_n$ in the basis $B = [e_1, \dots, e_n]$ which will be useful in the proof of Theorem \ref{theoremsmoothing}.

\begin{general}{Lemma}\label{lemmacoprimeN}
Let $a_1, \dots, a_n \in \Lambda$ be linearly independent. 
Let $\alpha_1, \dots, \alpha_{n}$ be the primitive positive dual family to $a_1, \dots, a_n$ in $L$. For all $1 \leq j \leq n$, write 
$$ \alpha_j = \sum_{j = 1}^{n} \alpha_{k, j} e_k$$
with $\alpha_{k,j} \in \zz$. The order of $\alpha_j \mod L'$ in $L/L'$ is precisely $n_j = N/\pregcd(N, \alpha_{1,j})$. In particular, the smoothing lattice $L'$ is \good for $a_1, \dots, a_n$ if and only if $\pregcd(N, \alpha_{1, j}) = 1$ for all $1 \leq j \leq n$.
\end{general}

\begin{proof}
Let $\alpha'_1, \dots, \alpha'_{n}$ be the primitive positive dual family to $a'_1, \dots, a'_{n}$ in $L'$. Then $\alpha'_1, \dots, \alpha'_{n}$ is also a positive dual family to $a_1, \dots, a_{n}$ in $L$. It follows from lemma [\hspace{1sp}\cite{firstpaper}, lemma 6] that there are positive integers $m_1, \dots, m_{n}$ such that $\alpha'_j = m_j \alpha_j$ for all $1 \leq j \leq n$ where $ \alpha_1, \dots, \alpha_{n}$ is the primitive positive dual family to $a_1, \dots, a_{n}$ in $L$. Write for all $1 \leq j \leq n$:
$$ \alpha_j = \sum_{k = 1}^{n} \alpha_{k, j} e_k$$
with $\alpha_{k,j} \in \zz$ and $\pregcd(\alpha_{1,j}, \dots, \alpha_{n,j}) = 1$. By definition of the integer $n_j = N/\pregcd(N, \alpha_{1,j}) = N/l_j$, the vector 
$$n_j\alpha_j = \frac{\alpha_{1,j}}{l_j} (Ne_1) + \sum_{k = 2}^{n} \frac{N\alpha_{k, j}}{l_j} e_k$$
belongs to $L'$ and $\pregcd(\alpha_{1,j}/l_j, N\alpha_{2, j}/l_j, \dots, N\alpha_{n, j}/l_j) = 1$. In particular, $n_j\alpha_j$ is a primitive vector in $L'$ and $\alpha'_j = (m_j/n_j) (n_j\alpha_j)$ is also a primitive vector in $L'$. Therefore $m_j = n_j$ and $\alpha'_j = n_j\alpha_j$. 
\end{proof}

For the rest of section \ref{sectionsmoothing} we suppose that the smoothing lattice $L'$ is indeed \good for the linear forms $a_1, \dots, a_n$ and we focus on the case where $a_1, \dots, a_n$ are linearly independent.
The goal of this section is to give an explicit formulation for the rational function $B_{n, a_1, \dots, a_n}(v)(w,x,L,L')$ in terms of periodic Bernoulli polynomials. Indeed, let us recall the definition of the classic Bernoulli polynomials using the generating series:
\begin{equation}\label{defpolbern}
\frac{e^{Xz}}{e^z-1} = \sum_{k \geq 0} B_k(X)\frac{z^{k-1}}{k!}.
\end{equation}
We may introduce the periodic versions of the Bernoulli polynomials $b_k(x) = B_k(x - \lfloor x \rfloor)$ for $x \in \rr$. In this section we shall prove the following:

\begin{general}{Proposition}\label{propbernexplicit}
Assume that $a_1, \dots, a_n \in \Lambda_N$ are linearly independent and that the smoothing lattice $L'$ is \good for $a_1, \dots, a_n$. Let $\alpha_1, \dots, \alpha_n$ be the positive dual basis to $a_1, \dots, a_n$. Fix a set $\mathcal{F}$ of representatives for $L/M$ where $M = \oplus_{j = 1}^n \zz \alpha_j$. Then there are explicit integers $r_j(\delta) \in \zz$ for all $1 \leq j \leq n$ and all $\delta = \sum_{j = 1}^n \delta_j\alpha_j/s_j \in \mathcal{F}$ such that:
$$B_{n, a_1, \dots, a_n}(v)(w,x,L, L')  = \epsilon \sum_{m = 0}^{n} \frac{w^m}{m!} \sum_{k_1 + \dots + k_{n} = n-m} \sum_{\delta \in \mathcal{F}} Y(k_1, \dots, k_{n}, v, \delta)\prod_{j = 1}^{n}\frac{x(\alpha_j)^{k_j-1}}{k_j!}$$
where $\epsilon = \signdet(a_1,\dots, a_n)$, $v = \sum_{j = 1}^n v_j\alpha_j/s_j$,
\begin{equation}\label{definitionY}
Y(k_1, \dots, k_{n}, v, \delta) = N\sum_{q \in Q} \prod_{j = 1}^{n} b_{k_j}\left(\frac{v_j + \delta_j + (r_j(\delta) +q_j)s_j}{Ns_j}\right)N^{k_j-1}- \prod_{j = 1}^{n} b_{k_j}\left(\frac{v_j+\delta_j}{s_j}\right)
\end{equation}
and
\begin{equation}\label{definitionQ}
Q = \left\{ (q_1, \dots, q_n) \in (\zsz{N})^{n} \setseparator \sum_{j = 1}^n q_j \alpha_{1, j} \equiv 0 \mod N  \right\}.
\end{equation}
\end{general}

To prove Proposition \ref{propbernexplicit} we shall first give an explicit description of both $B_{n, a_1, \dots, a_n}(v)(w,x,L)$ and $B_{n, a_1, \dots, a_n}(v)(w, x, L')$ individually in terms of periodic Bernoulli polynomials. Then, we shall show how to express any set of representatives $\mathcal{F}'$ of $L'/M'$ in terms of a fixed set of representatives $\mathcal{F}$ of $L/M$ and in terms of the set $Q$ given by \refp{definitionQ}. Let us start with the explicit description in terms of periodic Bernoulli polynomials:

\begin{general}{Lemma}\label{lemmabernoulliexplicit}
Fix $\mathcal{F}$ (resp. $\mathcal{F'}$) a set of representatives for $L/M$ (resp. $L'/M'$). Then:
\begin{align*}
B_{n, a_1, \dots, a_n}(v)(w,x,L) &= \epsilon \sum_{m = 0}^{n} \frac{w^m}{m!} \sum_{k_1 + \dots + k_{n} = n-m} \sum_{\delta \in \mathcal{F}} \prod_{j = 1}^{n} \frac{b_{k_j}(\frac{v_j + \delta_j}{s_j})x(\alpha_j)^{k_j-1}}{k_j!} \\
B_{n, a_1, \dots, a_n}(v)(w,x,L') &= \epsilon \sum_{m = 0}^{n} \frac{w^m}{m!} \sum_{k_1 + \dots + k_{n} = n-m} \sum_{\delta' \in \mathcal{F'}}\prod_{j = 1}^{n} \frac{b_{k_j}(\frac{v_j + \delta'_j}{Ns_j})N^{k_j-1}x(\alpha_j)^{k_j-1}}{k_j!}
\end{align*}
where $\epsilon = \signdet(a_1, \dots, a_n)$, $v = \sum_{j = 1}^n v_j \alpha_j/s_j$ and the sums range over integers $k_1 \geq 0, \dots, k_n \geq 0$.
\end{general}

\begin{proof}
We recall that by definition (see formula \refp{defbernoulliexplicit}) : 
\begin{align*}
B_{n, a_1, \dots, a_n}(v)(w,x, L) & = \epsilon \times \coefficient{\sum_{\delta \in (v+ L) \cap P(\abar)}\frac{e^{wt}e^{x(\delta)t}}{\prod_{j = 1}^n (1-e^{x(\alpha_j)t})}}{t}{0}\\
B_{n, a_1, \dots, a_n}(v)(w,x, L') & = \epsilon \times \coefficient{\sum_{\delta' \in (v+ L') \cap N.P(\abar)}\frac{e^{wt}e^{x(\delta')t}}{\prod_{j = 1}^n (1-e^{x(N\alpha_j)t})}}{t}{0}
\end{align*}
where $\epsilon = \signdet(a_1, \dots, a_n)$. Let us define for simplicity $F = \{\delta \in L \setseparator, v+\delta \in P(\abar)\}$ and $F' = \{ \delta' \in L' \setseparator, v+\delta' \in N.P(\abar)\}$. It is clear that the set $F$ (resp. $F'$) is the unique set of representatives for $L/M$ (resp. $L'/M'$) such that $v + F \subset P(\abar)$ (resp. $v + F' \subset N.P(\abar)$). Let us write these sets explicitly as:
\begin{align*}
F &= \left\{\sum_{j = 1}^n \frac{\delta_j\alpha_j}{s_j} \in L \setseparator \forall\, 1 \leq j \leq n, 0 \leq v_j + \delta_j < s_j, \delta_j \in \zz \right\}\\
F' &= \left\{\sum_{j = 1}^n \frac{\delta'_j.N\alpha_j}{Ns_j} \in L' \setseparator \forall\, 1 \leq j \leq n, 0 \leq v_j + \delta'_j < Ns_j, \delta'_j \in \zz \right\}
\end{align*}
so that using the definition of the classic Bernoulli polynomials (see \refp{defpolbern}) we get:
\begin{align*}
B_{n, a_1, \dots, a_n}(v)(w,x,L) &= \epsilon \sum_{m = 0}^{n} \frac{w^m}{m!} \sum_{k_1 + \dots + k_{n} = n-m} \sum_{\delta \in F} \prod_{j = 1}^{n} \frac{B_{k_j}(\frac{v_j + \delta_j}{s_j})x(\alpha_j)^{k_j-1}}{k_j!} \\
B_{n, a_1, \dots, a_n}(v)(w,x,L') &= \epsilon \sum_{m = 0}^{n} \frac{w^m}{m!} \sum_{k_1 + \dots + k_{n} = n-m} \sum_{\delta' \in F'} \prod_{j = 1}^{n} \frac{B_{k_j}(\frac{v_j + \delta'_j}{Ns_j})N^{k_j-1}x(\alpha_j)^{k_j-1}}{k_j!} 
\end{align*}
Since for all $\delta \in F$ (resp. $\delta' \in F'$) and all $1 \leq j \leq n$, $0 \leq (v_j + \delta_j)/s_j < 1$ (resp. $0 \leq (v_j + \delta_j)/(Ns_j) < 1$) we may rewrite this using the periodic Bernoulli polynomials and then replace the sets $F$ and $F'$ with any sets of representatives $\mathcal{F}$ and $\mathcal{F}'$ for $L/M$ and $L'/M'$ respectively. This gives:
\begin{align*}
B_{n, a_1, \dots, a_n}(v)(w,x,L) &= \epsilon \sum_{m = 0}^{n} \frac{w^m}{m!} \sum_{k_1 + \dots + k_{n} = n-m} \sum_{\delta \in \mathcal{F}} \prod_{j = 1}^{n} \frac{b_{k_j}(\frac{v_j + \delta_j}{s_j})x(\alpha_j)^{k_j-1}}{k_j!} \\
B_{n, a_1, \dots, a_n}(v)(w,x,L') &= \epsilon \sum_{m = 0}^{n} \frac{w^m}{m!} \sum_{k_1 + \dots + k_{n} = n-m} \sum_{\delta' \in \mathcal{F'}}\prod_{j = 1}^{n} \frac{b_{k_j}(\frac{v_j + \delta'_j}{Ns_j})N^{k_j-1}x(\alpha_j)^{k_j-1}}{k_j!}.
\end{align*}
Indeed, if $\tilde{\delta} = \sum_{j = 1}^n \tilde{\delta}_j\alpha_j/s_j \in L$ represents the same class in $L/M$ as $\delta$ then there are integers $m_1, \dots, m_n$ such that $\tilde{\delta} = \delta + \sum_{j = 1}^n m_j \alpha_j$. Therefore: 
$$b_{k_j}\left(\frac{v_j + \tilde{\delta}_j}{s_j}\right) = b_{k_j}\left(\frac{v_j + \delta_j}{s_j} + m_j \right) =  b_{k_j}\left(\frac{v_j + \delta_j}{s_j}\right)$$
and the explicit description above do not depend on the choice of representatives for $L/M$ and $L'/M'$. This completes the proof.
\end{proof}

Let us now give an explicit link between the quotient sets $L/M$ and $L'/M'$. Fix $\mathcal{F}$ and $\mathcal{F}'$ sets of representatives for $L/M$ and $L'/M'$. Let us remark that the identifications $\mathcal{F} \simeq L/M$ and $\mathcal{F}' \simeq L'/M'$ induce group structures on both $\mathcal{F}$ and $\mathcal{F}'$ defined respectively by:
$$\delta_1 * \delta_2 \equiv \delta_1 + \delta_2 \mod M, ~~ \delta'_1 * \delta'_2 \equiv \delta'_1 + \delta'_ 2 \mod M' $$ 
The neutral elements of $\mathcal{F}$ and $\mathcal{F}'$ correspond to the representative for the trivial classes $M$ and $M'$ in $L/M$ and $L'/M'$. We now relate $\mathcal{F}$ and $\mathcal{F}'$ for a \good smoothing lattice $L'$.

\begin{general}{Lemma}\label{lemmareprsnake}
Suppose that the smoothing lattice $L'$ is good for the linear forms $a_1, \dots, a_n$. Consider the map
$$ f := \begin{cases} \mathcal{F}' & \to \mathcal{F} \\
\sum_{j = 1}^n \frac{\delta'_j\alpha_j}{s_j} & \to \sum_{j = 1}^n \frac{\delta_j\alpha_j}{s_j} \mod M
\end{cases}$$ 
where the $\delta_j$'s are given by Euclidian division as $\delta'_j = q_j s_j+ \delta_j$ with $0 \leq v_j + \delta_j < s_j$ and $\delta_j \in \zz$. Here it is understood that $f(\delta')$ is the representative in $\mathcal{F}$ for the class $\sum_{j = 1}^n \delta_j\alpha_j/s_j \mod M$. The map $f$ is a $N^{n-1}$ to $1$ surjective group morphism and its kernel is isomorphic to the group:
$$Q = \left\{ (q_1, \dots, q_n) \in \zsz{N}^{n} \setseparator \sum_{j = 1}^n q_j \alpha_{1, j} \equiv 0 \mod N \right\}.$$
\end{general}

\begin{proof}
Let us first prove that the map $f$ is surjective. Consider $\delta = \sum_{j = 1}^n \frac{\delta_j\alpha_j}{s_j} \in \mathcal{F}$ and write in coordinates:
$$\delta = \sum_{j = 1}^n\frac{\delta_j}{s_j}\alpha_{1,j} e_1 + \sum_{k = 2}^n \sum_{j = 1}^n\frac{\delta_j}{s_j}\alpha_{k,j}e_k $$
The assumption that $\delta \in L$ is equivalent to $\sum_{j = 1}^n\frac{\delta_j}{s_j}\alpha_{k,j} \in \zz$ for all $1 \leq k \leq n$. We now wish to find integers $r_1, \dots, r_n \in \{0, \dots, N-1\}^n$ such that 
$$ \sum_{j = 1}^n\frac{\delta_j}{s_j}\alpha_{1,j} + \sum_{j = 1}^n r_j \alpha_{1,j} \in N\zz.$$
It follows from lemma \ref{lemmacoprimeN} that the integer $\alpha_{1,1}$ is coprime to $N$ so there is an integer $\beta \in \zz$ such that $\beta \alpha_{1,1} \equiv 1 \mod N$. Fix $r_1$ the remainder in the Euclidian division of the integer $-\beta\sum_{j = 1}^n\frac{\delta_j}{s_j}\alpha_{1,j}$ by $N$ satisfying $0 \leq r_1 < N$. Fix $r_2 = \dots = r_n = 0$. Then:
$$ \sum_{j = 1}^n\frac{\delta_j}{s_j}\alpha_{1,j} + \sum_{j = 1}^n r_j \alpha_{1,j} \equiv \left(\sum_{j = 1}^n\frac{\delta_j}{s_j}\alpha_{1,j}\right)(1-\beta\alpha_{1,1})  \equiv 0 \mod N$$
This shows that $\tilde{\delta} = \delta + \sum_{j = 1}^n r_j\alpha_j \in L'$. Let $\delta'$ be the representative in $\mathcal{F'}$ for the class $[\tilde{\delta}]$, that is $\delta' = \delta + \sum_{j = 1}^n (r_j + m_j N)\alpha_j$ for some integers $m_1, \dots, m_n$. Then it is clear that $f(\delta') = \delta$. Thus $f$ is surjective. In addition, the preimage $f^{-1}(\delta)$ is explicitly given by:
$$f^{-1}(\delta) = \left\{\left[\delta + \sum_{j = 1}^n r_j\alpha_j + \sum_{j = 1}^n q_j \alpha_j\right] \mod M' \setseparator (q_1, \dots, q_n)\in Q \right\}.$$
Indeed, if $f(\delta') = \delta$ then for all $1 \leq j \leq n$, $\delta'_j = (r_j+q_j) s_j + \delta_j$ for some $q_j \in \zz$. The condition that $\delta'_j \in L'$ is equivalent to $\sum_{j = 1}^n (\delta_j +s_j(r_j+q_j))\alpha_{1,j} \in N\zz$ which by definition of the $r_j$'s is equivalent to $\sum_{j = 1}^n q_j \alpha_{1,j} \in N\zz$ i.e. $(q_1, \dots, q_n) \in Q$. In particular, $\ker f \simeq Q$.

Let us now prove that $f$ is a group morphism. Consider $\delta'_1, \delta'_2 \in \mathcal{F}'$, define $\delta'_3 = \delta'_1 * \delta'_2$ and write for $k = 1, 2, 3$:
$$\delta'_k =  \sum_{j = 1}^n \frac{\delta'_{k,j}\alpha_j}{s_j} $$
The definition of the group law gives that
\begin{equation}\label{deltatrois}
\delta'_{3,j} = (\delta'_{1,j} + \delta'_{2,j}) + Ns_j m'_j \text{ for some } m'_j \in \zz 
\end{equation}
Define for $k = 1, 2, 3$ and $1 \leq j \leq n$ the unique integers $\delta_{k,j}$ satisfying:
$$0 \leq \delta_{k,j} < a_j(\alpha_j), ~~ \delta'_{k,j} = q_{k,j}s_j + \delta_{k,j}, q_{k, j} \in \zz $$
so that $f(\delta'_k) = \sum_{j = 1}^n \frac{\delta_{k,j}\alpha_j}{s_j}$ for $k = 1, 2, 3$. The definition of the group law on $\mathcal{F}$ gives:
$$f(\delta'_1) * f(\delta'_2) \equiv \sum_{j = 1}^n \frac{\delta_{1,j} + \delta_{2,j}\alpha_j}{s_j} \mod M$$
whereas $f(\delta'_3) \equiv \sum_{j = 1}^n \frac{\delta_{3,j}\alpha_j}{s_j} \mod M$. We shall prove that $\delta_{1,j}+\delta_{2,j} - \delta_{3,j} \in s_j\zz$ for $1 \leq j \leq n$. Indeed, by definition of $\delta_{k,j}$ we get:
$$\delta_{1,j} + \delta_{2, j} = \delta'_{1,j} - q_{1,j}s_j + \delta'_{2,j} - q_{2,j} s_j$$
Using \refp{deltatrois} we get:
$$ \delta_{1,j} + \delta_{2, j} = \delta'_{3,j} +Ns_jm'_j - q_{1,j}s_j - q_{2,j} s_j = \delta_{3,j} + s_j(q_{3,j} + Nm'_j - q_{1,j} - q_{2,j}) \equiv \delta_{3,j} \mod s_j\zz$$
This show that $\delta_1 + \delta_2 \equiv \delta_3 \mod M$ and so that:
$$f(\delta'_1 * \delta'_2) = f(\delta'_3) = f(\delta'_1)*f(\delta'_2).$$
Thus $f$ is a group morphism. Next, we prove that $f$ is an $N^{n-1}$ to $1$ map. This is given by the snake lemma for the following commutative diagram: 

\begin{center}
\begin{tabular}{C C C C C C C C C}
& & 0 & & 0 & & \ker f & & \\
& & \downarrow & & \downarrow & & \downarrow & & \\
0 &\rightarrow & M' & \rightarrow & L' & \rightarrow & \mathcal{F}' & \rightarrow &0 \\
& & \downarrow & & \downarrow & & \downarrow & & \\
0 &\rightarrow & M & \rightarrow & L & \rightarrow & \mathcal{F} & \rightarrow & 0 \\
& & \downarrow & & \downarrow & & \downarrow& & \\
& & (\zsz{N})^n & & \zsz{N} & & 0& & \\
\end{tabular}
\end{center}
\bigskip

\noindent The snake lemma gives a connecting group morphism $\ker f \to (\zsz{N})^n$ such that the following sequence is exact: 
$$0 \to 0 \to \ker f \to (\zsz{N})^n \to \zsz{N} \to 0 $$
As a consequence, $\cardinalshort{\ker f} = N^{n-1}$ and $f$ is an $N^{n-1}$ to $1$ map. 
\end{proof}

As an immediate consequence of this lemma we obtain the following corollary:

\begin{general}{Corollary}\label{corfprime}
For any set of representatives $\mathcal{F}$ for $L/M$, and any set of representatives $\mathcal{Q}$ for $Q$, the set
$$\mathcal{F}' := \{ \delta + \sum_{j = 1}^n (r_j(\delta) + q_j) \alpha_j \setseparator (q_1, \dots, q_n) \in \mathcal{Q} \} $$
is a set of representatives for $L'/M'$, where the $r_j(\delta)$'s are defined in the proof of lemma \ref{lemmareprsnake}.
\end{general}

We may now prove Proposition \ref{propbernexplicit} with this particular choice of sets of representatives for both $L/M$ and $L'/M'$. \bigskip

\begin{proofbis}{Proof of Proposition \ref{propbernexplicit}}
Let us fix $\mathcal{F}$ and $\mathcal{Q}$ two sets of representatives for $L/M$ and $Q$ respectively, as well as $\mathcal{F}'$ the set of representatives for $L'/M'$ given by Corollary \ref{corfprime}. It follows from lemma \ref{lemmabernoulliexplicit} that:
$$B_{n, a_1, \dots, a_n}(v)(w,x,L') = \epsilon \sum_{m = 0}^{n} \frac{w^m}{m!} \sum_{k_1 + \dots + k_{n} = n-m} \sum_{\delta' \in \mathcal{F'}}\prod_{j = 1}^{n} \frac{b_{k_j}(\frac{v_j + \delta'_j}{Ns_j})N^{k_j-1}x(\alpha_j)^{k_j-1}}{k_j!}$$
where $\epsilon = \signdet(a_1, \dots, a_n)$. Using Corollary \ref{corfprime} we may rewrite this using a double sum on $\mathcal{F}$ and $\mathcal{Q}$ as:
$$B_{n, a_1, \dots, a_n}(v)(w,x,L') = \epsilon \sum_{m = 0}^{n} \frac{w^m}{m!} \sum_{k_1 + \dots + k_{n} = n-m} \sum_{\delta \in \mathcal{F}}\sum_{q \in \mathcal{Q}}\prod_{j = 1}^{n} \frac{b_{k_j}(\frac{v_j + \delta_j + (r_j(\delta) + q_j)s_j}{Ns_j})N^{k_j-1}x(\alpha_j)^{k_j-1}}{k_j!}$$
Note that this expression is independent from the choice of representative set $\mathcal{Q}$ for $Q$. Thus, by definition of the smoothed function $B_{n, a_1, \dots, a_n}(v)(w,x,L,L')$ (see \refp{defsmoothedbn}) and using once again lemma \ref{lemmabernoulliexplicit}:
$$B_{n, a_1, \dots, a_n}(v)(w,x,L, L')  = \epsilon \sum_{m = 0}^{n} \frac{w^m}{m!} \sum_{k_1 + \dots + k_{n} = n-m} \sum_{\delta \in \mathcal{F}} Y(k_1, \dots, k_{n}, v, \delta)\prod_{j = 1}^{n}\frac{x(\alpha_j)^{k_j-1}}{k_j!}$$
where for all $k_1 \geq 0, \dots, k_n \geq 0$ such that $\sum_{j = 1}^n k_j \leq n$ and for all $\delta \in \mathcal{F}$:
$$ Y(k_1, \dots, k_{n}, v, \delta) = N\sum_{q \in Q} \prod_{j = 1}^{n} b_{k_j}\left(\frac{v_j + \delta_j + (r_j(\delta) + q_j)s_j}{Ns_j}\right)N^{k_j-1}- \prod_{j = 1}^{n} b_{k_j}\left(\frac{\delta_j}{s_j}\right).$$
\end{proofbis}

\subsection{Smoothed Bernoulli rational functions are rational valued}\label{sectionrational}

The goal of this section is to prove the following rationality statement for the smoothed $B_{n, a_1, \dots, a_n}(v)(w, x,L,L')$ functions.

\begin{general}{Proposition}\label{propbernrational}
Assume that $a_1, \dots, a_n \in \Lambda_N$ are linearly independent and that the smoothing lattice $L'$ is \good for $a_1, \dots, a_n$. Let $\alpha_1, \dots, \alpha_n$ be the positive dual basis to $a_1, \dots, a_n$. Put as before $\epsilon = \signdet(a_1, \dots, a_n)$ and fix a set $\mathcal{F}$ of representatives for $L/M$. Then:
$$B_{n, a_1, \dots, a_n}(v)(w,x,L, L')  = \epsilon \sum_{\delta \in \mathcal{F}}\left(N\sum_{q \in Q} \prod_{j = 1}^n b_1\left(\frac{v_j + \delta_j + (r_j(\delta) +q_j)s_j}{Ns_j}\right) - \prod_{j= 1}^n b_1\left(\frac{v_j+\delta_j}{s_j}\right)\right)$$
where the integers $r_j(\delta)$ are given in the proof of lemma \ref{lemmareprsnake} and $Q$ is defined in Proposition \ref{propbernexplicit}. In particular $B_{n, a_1, \dots, a_n}(v)(w, x, L,L')$ is a rational number which depends only on the linear forms $a_1, \dots, a_n$ and on the class of $v$ in $V/L'$ but not on $w, x \in \cc \times \homlc$.
\end{general}

This Proposition essentially expresses the smoothed function $\smoothedbn$ in terms of 
a smoothed higher Dedekind sum. To prove this statement we shall prove that each term $Y(k_1, \dots, k_n, v, \delta)$ vanishes, unless $k_1 = k_2 = \dots = k_n = 1$. This is exactly the claim of the following crucial lemma:

\begin{general}{Lemma}\label{lemmayzero}
Suppose that $k_1, \dots, k_n$ are non-negative integers such that $\sum_{j = 1}^n k_j \leq n$. Suppose that $(k_1, \dots, k_n) \neq (1, \dots, 1)$. Then $Y(k_1, \dots, k_n, v, \delta) = 0$ for any $\delta \in \mathcal{F}$.
\end{general}

\begin{proof}
Let us fix non-negative integers $k_1, \dots, k_n$ such that $\sum_{j = 1}^n k_j \leq n$. The condition $(k_1, \dots, k_n) \neq (1, \dots, 1)$ is equivalent to the existence of an index $1 \leq j \leq n$ such that $k_j = 0$. Define $J = \{1 \leq j \leq n \setseparator k_j \neq 0 \}$ and $J^c = \{1, \dots, n\} - J \neq \emptyset$. Let us define a map:
$$g_{J} := \begin{cases} Q &\to \prod_{j \in J} \zsz{N} = Q(J) \\ (q_1, \dots, q_n) & \to (q_j)_{j \in J} \end{cases} $$
This map is clearly a group morphism, and we shall prove that it is surjective. Indeed, by assumption, $J^c$ is not empty so we may fix an index $j' \in J^c$. The condition that the smoothing lattice $L'$ is good for the linear forms $a_1, \dots, a_n$ implies that $\alpha_{1, j'}$ is coprime to $N$. Therefore, for any element $(q_j)_{j \in J} \in Q(J)$ there is an integer $q_{j'} \in \zsz{N}$ such that $\sum_{j \in J} q_j \alpha_{1,j} + q_{j'}\alpha_{1, j'} \equiv 0 \mod N$. Setting $q_{j''} = 0$ for all $j'' \in J^c - \{j'\}$ gives $(q_1, \dots, q_n) \in Q$ and $g_{J}(q_1, \dots, q_n) = (q_j)_{j \in J}$. Therefore, $g_{J}$ is surjective. We shall use the function $g_{J}$ to switch the sum and product in the expression of $Y(k_1, \dots, k_n, \delta)$. First, let us rewrite $Y(k_1, \dots, k_n, \delta)$ in terms of the set $J$:
\begin{align*}
Y(k_1, \dots, k_{n}, v, \delta) & = N^{1-n}\sum_{q \in Q} \prod_{j = 1}^{n} b_{k_j}\left(\frac{v_j+\delta_j+(r_j(\delta) + q_j)s_j}{Ns_j}\right)N^{k_j}- \prod_{j = 1}^{n} b_{k_j}\left(\frac{v_j+\delta_j}{s_j}\right) \\
Y(k_1, \dots, k_{n}, v, \delta) & = N^{1-n}\sum_{q \in Q} \prod_{j \in J} b_{k_j}\left(\frac{v_j+\delta_j+(r_j(\delta) + q_j)s_j}{Ns_j}\right)N^{k_j}- \prod_{j \in J} b_{k_j}\left(\frac{v_j+\delta_j}{s_j}\right)
\end{align*}
where we have used the fact that if $j \not\in J$ then $b_{k_j} = b_0$ is the constant function equal to $1$. From this expression it is clear that the term 
$$\prod_{j \in J} b_{k_j}\left(\frac{v_j+\delta_j+(r_j(\delta) + q_j)s_j}{Ns_j}\right)N^{k_j} $$
only depends on the image of $g(q) \in Q_J$ of $q$, therefore:
$$Y(k_1, \dots, k_{n}, v, \delta) = N^{1-n}\cardinalshort{\ker(g_{J})}\sum_{q_J \in Q(J)} \prod_{j \in J} b_{k_j}\left(\frac{v_j+\delta_j+(r_j(\delta) +q_j)s_j}{Ns_j}\right)N^{k_j}- \prod_{j \in J} b_{k_j}\left(\frac{v_j+\delta_j}{s_j}\right)$$
We may now use the fact that $\cardinalshort{\ker(g_j)} = N^{n-1-\cardinalshort{J}}$ and then the fact that $Q_J = \prod_{j \in J} \zsz{N}$ to switch sum and product in the expression of $Y(k_1, \dots, k_n, \delta)$ to obtain:
\begin{align*}
Y(k_1, \dots, k_{n}, v, \delta) & = \sum_{q_J \in Q(J)} \prod_{j \in J} b_{k_j}\left(\frac{v_j+\delta_j+(r_j(\delta) + q_j)s_j}{Ns_j}\right)N^{k_j-1}- \prod_{j \in J} b_{k_j}\left(\frac{v_j+\delta_j}{s_j}\right) \\
Y(k_1, \dots, k_{n}, v, \delta) & = \prod_{j \in J} \sum_{q_j \in \zsz{N}} b_{k_j}\left(\frac{v_j+\delta_j+(r_j(\delta)+q_j)s_j}{Ns_j}\right)N^{k_j-1}- \prod_{j \in J} b_{k_j}\left(\frac{v_j+\delta_j}{s_j}\right)
\end{align*}
It follows from the well-known distribution relation:
\begin{equation}\label{eqdistribbern}
\sum_{k \in \zsz{N}} b_m\left(\frac{x+k}{N}\right) N^{m-1}= b_m(x)
\end{equation}
applied here to $x = (v_j + \delta_j + r_j(\delta)s_j)/s_j$ that:
$$Y(k_1, \dots, k_{n}, v, \delta) = \prod_{j \in J} b_{k_j}\left(\frac{v_j+\delta_j + r_j(\delta)s_j}{s_j}\right)- \prod_{j \in J} b_{k_j}\left(\frac{v_j+\delta_j}{s_j}\right).$$
Since the $b_{k_j}$ functions are $1$-periodic and $r_j \in \zz$ we get the desired conclusion 
$$Y(k_1, \dots, k_n, v, \delta) = 0.$$
\end{proof}

We are now ready to deduce Proposition \ref{propbernrational} from lemma \ref{lemmayzero}. \bigskip

\begin{proofbis}{Proof of Proposition \ref{propbernrational}}
It follows from Proposition \ref{propbernexplicit} that:
$$B_{n, a_1, \dots, a_n}(v)(w,x,L, L')  = \epsilon\sum_{m = 0}^{n} \frac{w^m}{m!} \sum_{k_1 + \dots + k_{n} = n-m} \sum_{\delta \in \mathcal{F}} Y(k_1, \dots, k_{n}, v, \delta)\prod_{j = 1}^{n}\frac{x(\alpha_j)^{k_j-1}}{k_j!}$$
where $\epsilon = \signdet(a_1, \dots, a_n)$. Since for all $(k_1, \dots, k_n) \neq (1, \dots, 1)$ and all $\delta \in \mathcal{F}$, $Y(k_1, \dots, k_n, v, \delta) = 0$ the sum reduces to: 
$$B_{n, a_1, \dots, a_n}(v)(w,x,L, L') = \epsilon\sum_{\delta \in \mathcal{F}} Y(1, \dots, 1, v, \delta)$$
which gives exactly the desired expression when replacing $Y(1, \dots, 1, v, \delta)$ by its definition:
$$B_{n, a_1, \dots, a_n}(v)(w,x,L, L')  = \epsilon\sum_{\delta \in \mathcal{F}}\left(N\sum_{q \in Q} \prod_{j = 1}^n b_1\left(\frac{v_j + \delta_j + (r_j(\delta) +q_j)s_j}{Ns_j}\right) - \prod_{j= 1}^n b_1\left(\frac{v_j+\delta_j}{s_j}\right)\right).$$
It is then clear that $B_{n, a_1, \dots, a_n}(v)(w,x,L, L')$ is a rational number depending only on the linear forms $a_1, \dots, a_n$ and on the class of $v$ in $V/L'$ but not on $w, x \in \cc \times \homlc$.
\end{proofbis}

\subsection{Smoothed Bernoulli rational functions have bounded denominators}\label{sectiontraces}

In this section we carry out the proof of Theorem \ref{theoremsmoothing} by expressing the smoothed functions $B_{n, a_1, \dots, a_n}(v)(w,x,L,L')$ in terms of traces of cyclotomic units, thus proving a uniform bound on their denominator in terms of the dimension $n$ and the smoothing index $N$. To achieve this, we will use the Fourier transformation on the finite group $\zsz{N}$ to rewrite the expression obtained in Proposition \ref{propbernrational} borrowing ideas from \cite{CD}. Let us denote by $\zeta = \zeta_N = \exp(2i\pi/N)$ a primitive $N$-th root of unity.
Introduce the auxiliary function $\chi : (\zsz{N})^{n} \to \{0, N\}$ defined by:
$$ \chi(q) = \sum_{k \in \zsz{N}} \zeta^{q_j\alpha_{1, j}k} = \begin{cases} 0 \text{ if } q \not\in Q \\ N \text{ if } q \in Q\end{cases}$$
where $\alpha_{1,j} = \langle \alpha_j, e_1 \rangle$. We may then write following Proposition \ref{propbernrational} and using the auxiliary function $\chi$:
\begin{equation}\label{eqfourierz}
B_{n, a_1, \dots, a_n}(v)(w,x,L, L') =  \epsilon\sum_{\delta \in \mathcal{F}} \left(Z(\delta)- \prod_{j = 1}^{n} b_1\left(\frac{v_j+\delta_j}{s_j}\right)\right)
\end{equation}
where 
$$Z(\delta) := \sum_{q \in \zsz{N}} \chi(q) \prod_{j = 1}^{n} b_1\left(\frac{v_j+\delta_j + (r_j(\delta) + q_j)s_j}{Ns_j}\right).$$
Using the definition of $\chi$ we may rewrite this auxiliary function $Z(\delta)$ as:
$$Z(\delta) = \sum_{q \in \zsz{N}^n}\sum_{k \in \zsz{N}}\prod_{j = 1}^{n} \zeta^{q_j\alpha_{1,j}k}b_1\left(\frac{v_j+\delta_j + (r_j(\delta)+q_j)s_j}{Ns_j}\right)$$
and we remark that the sum over $q \in \zsz{N}^n$ may be then be inverted with the product over $1 \leq j \leq n$ so:
\begin{equation}\label{eqfourierref}
Z(\delta) = \sum_{k \in \zsz{N}}\prod_{j = 1}^{n}\sum_{q_j \in \zsz{N}}\zeta^{q_j\alpha_{1,j}k}b_1\left(\frac{v_j+\delta_j + (r_j(\delta) + q_j)s_j}{Ns_j}\right).\end{equation}
The proof of Theorem \ref{theoremsmoothing} will essentially follow from the following lemma which is a reformulation of [\hspace{1sp}\cite{CD}, Lemma 2.13].

\begin{general}{Lemma}\label{lemmabernCD}
Suppose $N \geq 2$ is an integer. If $x \in \rr$ and $y \in \bb{F}_N -\{0\}$:
$$\sum_{q \in \zsz{N}} \zeta^{yq}b_1\left(\frac{x+q}{N}\right) = \frac{\zeta^{-y\lfloor x\rfloor}}{\zeta^y-1} $$ 
where $\zeta = \exp(2i\pi/N)$.
\end{general}

This lemma allows us to write the expression $Z(\delta)$ in terms of traces of cyclotomic units, and we may now prove Theorem \ref{theoremsmoothing}.

\begin{proofbis}{Proof of Theorem \ref{theoremsmoothing}}
Let us first treat the term 
$$\sum_{q_j \in \zsz{N}} \zeta^{q_j\alpha_{1,j}k}b_1\left(\frac{v_j+\delta_j + (r_j(\delta) + q_j)s_j}{Ns_j}\right)$$ 
in expression \refp{eqfourierref} for $k = 0$ and $1 \leq j \leq n$ using the distribution relation \refp{eqdistribbern}. This gives:
$$\sum_{q_j \in \zsz{N}} b_1\left(\frac{v_j+\delta_j + (r_j(\delta) + q_j)s_j}{Ns_j}\right) = N^{0}b_1\left(\frac{v_j+\delta_j + r_j(\delta)s_j}{s_j}\right) = b_1\left(\frac{v_j+\delta_j}{s_j}\right)$$
since $r_j(\delta) \in \zz$ and $b_1$ is $1$-periodic. Thus the term for $k = 0$ cancels with the term $\prod_{j = 1}^n b_1\left(\frac{v_j+\delta_j}{s_j}\right)$ in expression \refp{eqfourierz} and:
\begin{equation}\label{eqfourierrefdeux}
B_{n, a_1, \dots, a_n}(v)(w,x,L, L') = \epsilon\sum_{\delta \in \mathcal{F}} \sum_{k = 1}^{N-1}\prod_{j = 1}^{n}\sum_{q_j \in \zsz{N}}\zeta^{q_j\alpha_{1,j}k}b_1\left(\frac{v_j+\delta_j + (r_j(\delta) + q_j)s_j}{Ns_j}\right)
\end{equation}
Let us now apply lemma \ref{lemmabernCD} to each term in expression \refp{eqfourierrefdeux} to obtain:
$$\sum_{q_j \in \zsz{N}} \zeta^{q_j\alpha_{1,j}k}b_1\left(\frac{v_j+\delta_j + (r_j(\delta) + q_j)s_j}{Ns_j}\right) =  \left(\frac{\zeta^{-\alpha_{1,j}k\lfloor\frac{v_j + \delta_j + r_j(\delta)s_j}{s_j}\rfloor}}{\zeta^{\alpha_{1,j}k}-1}\right). $$ 
Since $r_j(\delta) \in \zz$ we get $\lfloor\frac{v_j + \delta_j + r_j(\delta)s_j}{s_j}\rfloor = \lfloor\frac{v_j + \delta_j}{s_j}\rfloor$ and therefore:
$$ B_{n, a_1, \dots, a_n}(v)(w,x,L, L') =  \epsilon\sum_{\delta \in \mathcal{F}} \left(\sum_{k = 1}^{N-1}\prod_{j = 1}^{n} \left(\frac{\zeta^{-\alpha_{1,j}k\lfloor\frac{v_j + \delta_j}{s_j}\rfloor}}{\zeta^{\alpha_{1,j}k}-1}\right)\right)$$
This may be written as a sum of traces of cyclotomic units using the well-known bijection:
$$\zsz{N} - \{0\} = \disjointunion_{d | N, d \neq 1} \zsz{d}^{\times} $$
which gives:
$$ B_{n, a_1, \dots, a_n}(v)(w,x,L, L') =  \epsilon\sum_{\delta \in \mathcal{F}}\sum_{d | N, d \neq 1} \trace_{\qq(\zeta_d)/\qq}\left(\prod_{j = 1}^{n} \left(\frac{\zeta_d^{-\alpha_{1,j}\lfloor\frac{v_j + \delta_j}{s_j}\rfloor}}{\zeta_d^{\alpha_{1,j}}-1}\right)\right)$$
where for all $d|N$, $\zeta_d = \exp(2i\pi/d)$ and $\trace_{\qq(\zeta_d)/\qq}$ is the trace from $\qq(\zeta_d)$ to $\qq$. This is the desired relation and this completes the proof of Theorem \ref{theoremsmoothing}.
\end{proofbis}

We end this section by using Theorem \ref{theoremsmoothing} to obtain a uniform bound on the denominators of all $B_{n, a_1, \dots, a_n}(v)(w,x,L,L')$ functions in terms of $n$ and $N$, thus proving Theorem \ref{maintheorem}. Results of this type are now classic and some may be found in \cite{ZagierDedekind}, \cite{DasguptaShintani} and \cite{CD}. \medskip

\begin{proofbis}{Proof of Theorem \ref{maintheorem}}
We first make the remark that when $a_1, \dots, a_n \in \Lambda_N$ are linearly dependent and in good position in $\dualsimple{V}$, the function $B_{n, a_1, \dots, a_n}$ is identically $0$ and we may set $b(a_1, \dots, a_n, v) = 0$ identically in that case. Let us now suppose that $a_1, \dots, a_n \in \Lambda_N$ are linearly independent and that the smoothing lattice $L'$ is \good for $a_1, \dots, a_n$. Let us study each term $\trace_{\qq(\zeta_d)/\qq}(u_d)$ appearing in Theorem \ref{theoremsmoothing}, where 
$$ u_d = \prod_{j = 1}^{n} \left(\frac{\zeta_d^{-\alpha_{1,j}\lfloor\frac{v_j + \delta_j}{s_j}\rfloor}}{\zeta_d^{\alpha_{1,j}}-1}\right)$$
for $d | N, d \neq 1$. On the one hand, if $d$ is divisible by two distinct primes, as $\alpha_{1,j}$ is coprime to $d$, it is well-known that $\zeta_d^{\alpha_{1,j}}-1$ is a unit inside $\qq(\zeta_d)$ which implies that $u_d$ is a unit in the ring of integers $\mathcal{O}_{\qq(\zeta_d)}$ of $\qq(\zeta_d)$ and $\trace_{\qq(\zeta_d)/\qq}(u_d) \in \zz$. On the other hand, when $d = p^{\nu}$ is a power of the prime $p$, the algebraic integer $\zeta_d^{\alpha_{1,j}}-1$ is a generator for the unique prime ideal $\goth{P}_d$ above $p$ in $\qq(\zeta_d)$. This cyclotomic extension is totally ramified at $p$, therefore $\goth{P}_d^{\varphi(d)} = (p)$ where $\varphi(d) = p^{\nu-1}(p-1)$ is Euler's totient function evaluated at $d = p^{\nu}$. Thus $u_d$ is a generator of the ideal $\goth{P}_d^{-n}$. Let us introduce the different ideal $\mathfrak{D}$ of $\qq(\zeta_d)$ defined by:
$$\mathfrak{D}^{-1} = \{ x \in \qq(\zeta_d) \setseparator \forall y \in \mathcal{O}_{\qq(\zeta_d)}, \trace_{\qq(\zeta_d)/\qq}(xy) \in \zz\}.$$
In particular, any element $u$ in $\mathfrak{D}^{-1}$ satisfies $\trace_{\qq(\zeta_d)/\qq}(u) = \trace_{\qq(\zeta_d)/\qq}(u \times 1) \in \zz$.
It follows from [\hspace{1sp}\cite{Neukirch}, Lemma 10.1] that in this particular situtation $\goth{P}_d^{m}$ divides exactly $\mathfrak{D}$ where $m = p^{\nu-1}(p\nu -\nu-1)$. If we find an integer $k$ such that $p^k u_d \in \mathfrak{D}^{-1}$ then we will obtain $\trace_{\qq(\zeta_d)/\qq}(p^k u_d) \in \zz$ and therefore $\trace_{\qq(\zeta_d)/\qq}(u_d) \in p^{-k}\zz$. To determine the minimal such integer $k$ let us remark that:
$$ (p^k u_d) = \goth{P}^{k\varphi(p^{\nu}) -n} $$
therefore $p^k u_d \in \mathfrak{D}^{-1}$ if and only if $k \varphi(p^{\nu}) -n \geq -p^{\nu-1}(p\nu - \nu -1)$ which gives the condition:
$$k \geq \frac{n}{\varphi(p^{\nu})} - 1 + \frac{1}{p-1}.$$
Therefore:
$$\trace_{\qq(\zeta_d)/\qq}(u_d) \in p^{- \lceil\frac{n}{\varphi(d)}-1+\frac{1}{p-1}\rceil}\zz$$
where $\lceil x \rceil$ is the ceiling function satisfying $\lceil x \rceil -1 < x \leq \lceil x \rceil$ and $\lceil x \rceil \in \zz$. We can then conclude that 
$$B_{n, a_1, \dots, a_n}(v)(w,x,L, L') \in \sum_{p | N} \sum_{\nu = 1}^{v_p(N)} p^{-\left\lceil\frac{n}{(p-1)p^{\nu-1}}-1 + \frac{1}{p-1} \right\rceil} \zz$$
where $v_p(N)$ is the $p$-adic valuation of $N$. For all prime divisor $p$ of $N$ the term for $\nu = 1$ is dominant therefore we get the simpler relation:
$$B_{n, a_1, \dots, a_n}(v)(w,x,L, L') \in \sum_{p | N} p^{-\left\lceil\frac{n+1}{p-1} -1 \right\rceil} \zz.$$
It is not hard to check that $\left\lceil\frac{n+1}{p-1} -1 \right\rceil = \lfloor \frac{n}{p-1} \rfloor$ for any $n \geq 2, p \geq 2$, therefore we may set $\denomclasscohom = \prod_{p | N} p^{\left\lfloor\frac{n}{p-1} \right\rfloor}$, which uniformly bounds the denominators of all values $B_{n, a_1, \dots, a_n}(v)(w, x, L,L')$:
$$B_{n, a_1, \dots, a_n}(v)(w,x,L, L') \in \denomclasscohom^{-1} \zz.$$
Lastly, define $b(a_1, \dots, a_n,v)$ to be precisely the integer $B_{n, a_1, \dots, a_n}(v)(w,x,L,L').\denomclasscohom$. It follows from Proposition \ref{propbernrational} that $b(a_1, \dots, a_n,v)$ does not depend on the choice of $w, x \in \cc \times \homlc$. Putting this together with Theorem \ref{theoremmodularextend} shows that when $a_1, \dots, a_n$ are \wellplaced in $\dualsimple{V}$ and the smoothing lattice $L'$ is \good for $a_1, \dots, a_n$: 
$$\left(\prod_{j = 1}^{n} G_{n-2, a_1, \dots, \omitvar{a_j}, \dots, a_n}(v)(w, \plgt, L, L')^{(-1)^{j+1}}\right)= \exp\left(\frac{2i\pi b(a_1, \dots, a_n, v)}{\denomclasscohom}\right) $$
which is the desired result.
\end{proofbis}

\section{Cohomological interpretation}\label{sectionsmoothingcocycle}

In this last section we give a cohomological interpretation of the results presented in this paper. We first recall some functions introduced in \cite{firstpaper} which satisfy cocycle and coboundary properties as a consequence of formulae \refp{modularproperty} and \refp{cocyclebernoulli}. Then, we introduce smoothed versions of these functions and show how the smoothing operation affects their properties. Lastly, we restrict these functions to tori associated to groups of units in number fields, yielding proper cocycles for subgroups of $\congruencegrouphnn$.

In the first article in this series, we defined two collections of functions $\psi_{n,a}$ and $\phi_{n,a}$ attached to a primitive linear form $a \in \Lambda$ defined by:
$$\psi_{n,a}  := \begin{cases}\slnz{n}^{n-2} & \to \mathcal{F}(V/L \times \cc \times \homlc, \cc) \\ (g_1, \dots, g_{n-2}) & \to \left((v, w, \plgt) \to G_{n-2, a, g_1\cdot a, \dots, (g_1\dots g_{n-2}) \cdot a}(v)(w, \plgt)\right)\end{cases}$$
$$\phi_{n,a}  := \begin{cases} \slnz{n}^{n-1} & \to \mathcal{F}(V/L, \qq[w](\plgt)) \\ (g_1, \dots, g_{n-1}) & \to  B_{n, a, g_1\cdot a, (g_1g_2)\cdot a, \dots, (g_1\dots g_{n-1})\cdot a}(v)(w, \plgt)\end{cases}.$$
The modular property \refp{modularproperty} was then rephrased by saying that when the linear forms $a, g_1 \cdot a, \dots, (g_1\dots g_{n-1})\cdot a$ were linearly independent, the multiplicative coboundary of $\psi_{n,a}$ was given by $\exp(2i\pi\phi_{n, a})$ as:
\begin{equation}\label{eqcohomdpsi}\cohomdx\psi_{n, a}(g_1, \dots, g_{n-1}) = \exp(2i\pi \phi_{n, a}(g_1, \dots, g_{n-1})).\end{equation}
Theorem \ref{theoremmodularextend} implies that the coboundary relation \refp{eqcohomdpsi} holds for any $a_1, \dots, a_n$ which are \wellplaced in $\dualsimple{V}$, that is whenever $\rank(a_1, \dots, a_n) \neq n-1$ or whenever $\rank(a_1, \dots, a_n) = n-1$ and $0$ is not a barycenter of $a_1, \dots, a_n$ in $\dualsimple{V}$. On the other hand, the cocycle relation \refp{cocyclebernoulli} gives the partial cocycle relation:
\begin{equation}\label{eqcohomdphi}\cohomd\phi_{n,a}(g_1, \dots, g_{n}) = 0 \end{equation}
for any $g_1, \dots, g_{n} \in \slnz{n}$ such that $a, g_1 \cdot a, \dots, (g_1\dots g_n) \cdot a$ are not in bad position in $\dualsimple{V}$ in the sense of \cite{firstpaper}. 

Under the smoothing operation we introduced in this article we may define smoothed versions of these functions for a primitive linear form $a \in \Lambda_N$:
$$\psi_{n,a}^{(N)}  := \begin{cases}\congruencegrouphnn^{n-2} & \to \mathcal{F}(V/L' \times \cc \times \homlc, \cc) \\ (g_1, \dots, g_{n-2}) & \to \left((v, w, \plgt) \to G_{n-2, a, g_1\cdot a, \dots, (g_1\dots g_{n-2}) \cdot a}(v)(w, \plgt, L, L')\right)\end{cases}$$
$$\phi_{n,a}^{(N)}  := \begin{cases} \congruencegrouphnn^{n-1} & \to \mathcal{F}(V/L', \qq[w](\plgt)) \\ (g_1, \dots, g_{n-1}) & \to  B_{n, a, g_1\cdot a, (g_1g_2)\cdot a, \dots, (g_1\dots g_{n-1})\cdot a}(v)(w, \plgt, L, L')\end{cases}.$$
Let us now rephrase the main results of this paper in terms of these smoothed functions. It follows from \refp{eqcohomdpsi} that they satisfy the coboundary relation
\begin{equation}\label{eqcohomdpsiN}\cohomdx\psi_{n, a}^{(N)}(g_1, \dots, g_{n-1}) = \exp(2i\pi \phi_{n, a}^{(N)}(g_1, \dots, g_{n-1}))\end{equation}
whenever $a, g_1 \cdot a, \dots, (g_1 \dots g_{n-1})\cdot a$ are \wellplaced in $\dualsimple{V}$. The cocycle relation \refp{eqcohomdphi} also directly gives the cocycle relation
\begin{equation}\label{eqcohomdphiN}\cohomd\phi_{n,a}^{(N)}(g_1, \dots, g_{n}) = 0 \end{equation}
for any $g_1, \dots, g_{n} \in \congruencegrouphnn$ such that $a, g_1 \cdot a, \dots, (g_1\dots g_n) \cdot a$ are not in bad position in $\dualsimple{V}$ in the sense of \cite{firstpaper}. In addition, Theorem \ref{theoremsmoothing} implies that whenever the smoothing lattice $L'$ is \good for the linear forms $a, g_1 \cdot a, \dots, (g_1\dots g_{n-1}) \cdot a$, the value of the function $\phi_{n,a}^{(N)}(g_1, \dots, g_{n-1})$ does not depend on $w, x \in \cc \times \homlc \simeq \cc \times \cc^n$ and thus the function $\phi_{n, a}^{(N)}(g_1, \dots, g_{n-1})$ may be viewed as an element of $\mathcal{F}(V/L', \denomclasscohom^{-1} \zz)$ which is essentially a smoothed higher Dedekind sum. As for the function $\psi_{n,a}^{(N)}$, it follows from Corollary \ref{cormaintheorem} that for any $g_1, \dots, g_{n-1} \in \congruencegrouphnn$:
$$\cohomdx \left(\psi_{n, a}^{(N)}\right)^{\denomclasscohom} = 1$$
whenever the linear forms $a, g_1\cdot a, \dots, (g_1 \dots g_{n-1})\cdot a$ satisfy the hypothesis of Theorem \ref{maintheorem}.

Let us now consider specific subgroups $U$ of $\congruencegrouphnn$ satisfying the following property: $\forall m \geq 2, \forall g_1, \dots, g_m \in U, \forall \mu_1, \dots, \mu_m \in \zz_{\geq 0}$,
\begin{equation}\label{conditionsubgroup}
\sum_{j  = 1}^m \mu_j(g_j \cdot a) = 0 \Rightarrow \mu_1 = 0, \dots, \mu_m = 0
\end{equation}
(see [\hspace{1sp}\cite{firstpaper}, Condition (26)]). In particular, this condition is satisfied by certain unit groups in number fields as explained in \cite{firstpaper}. Indeed, if $\kk$ is a number field of degree $n$ with at least one real place $\sigma_{\rr}$, we may consider subgroups $\mathcal{U}$ of the unit group $\units$ satisfying $\forall\,\eps \in \mathcal{U}, \sigma_{\rr}(\eps) > 0$. We may consider lattices $L$ and $L'$ in $\kk$ corresponding to fractional ideals of $\kk$ such that $L' \subset L$ and $L/L' \simeq \zsz{N}$. In particular, these fractional ideals are stabilised by elements of $\mathcal{U}$ so  we may identify each element $\eps \in \mathcal{U}$ with the matrix $M_{\eps}$ corresponding to the multiplication by $\eps$ in a basis $B' = [Ne_1, e_2, \dots, e_n]$ of $L'$ such that $B = [e_1, e_2, \dots, e_n]$ is a $\zz$-basis of $L$. The group $U = \{ M_{\eps} \setseparator \eps \in \mathcal{U}\} \simeq \mathcal{U}$ is then an abelian subgroup of $\congruencegrouphnn$ satisfying \refp{conditionsubgroup}. It follows from [\hspace{1sp}\cite{firstpaper}, \symbolparagraph 3.3] that the application $\phi_{n,a}^{(N)}$ is a $(n-1)$-cocycle on the group $U \simeq \mathcal{U}$ with values in $\mathcal{F}(V/L', \qq[w](\plgt))$ and Theorem \ref{theoremsmoothing} implies that  $\phi_{n,a}^{(N)}(g_1, \dots, g_{n-1}) \in \mathcal{F}(V/L', \denomclasscohom^{-1} \zz)$ whenever the smoothing lattice is \good for $a, g_1 \cdot a, \dots, (g_1, \dots, g_{n-1})\cdot a$. Under hypothesis \refp{conditionsubgroup} it is also true that the splitting relation \refp{eqcohomdpsiN} holds for all $g_1, \dots, g_{n-1} \in U$:
$$\cohomdx\psi_{n, a}^{(N)}(g_1, \dots, g_{n-1}) = \exp(2i\pi \phi_{n, a}^{(N)}(g_1, \dots, g_{n-1}))$$
as \refp{conditionsubgroup} guarantees that the linear forms $a, g_1 \cdot a, \dots, (g_1 \dots g_{n-1})\cdot a$ are \wellplaced in $\dualsimple{V}$.

It would be interesting to find some conditions on the linear form $a \in \Lambda_N$ and on the unit group $\mathcal{U}$ such that for a fixed smoothing lattice $L' \subset L$, the lattice $L'$ is \good for any family $g_1 \cdot a, \dots, g_n \cdot a$, where $g_1, \dots, g_n \in \mathcal{U}$ (i.e. the smoothing lattice is ``uniformly \good for $a$ and $\mathcal{U}$''). Indeed, if the smoothing lattice $L'$ were to be uniformly \good for $a$ and $\mathcal{U}$ then the restriction of the smoothed function $\phi_{n,a}^{(N)}$ to $U \simeq \mathcal{U}$ would give a cocycle in $H^{n-1}(\mathcal{U}, \mathcal{F}(V/L', \denomclasscohom^{-1}\zz))$. In addition, the restriction of $(\psi_{n,a}^{(N)})^{\denomclasscohom}$ to $U \simeq \mathcal{U}$ would yield a multiplicative cocycle in $H^{n-2}(\mathcal{U}, \mathcal{F}(V/L' \times \cc \times \homlc, \cc))$. We discuss this question of finding a condition ensuring that the smoothing lattice $L'$ is uniformly \good for $a$ and $\mathcal{U}$ in the third paper in this series, where we consider groups of totally positive units in $\ok$ which are congruent to $1$ modulo some integral ideal $\goth{f}$, where $\kk$ is a number field of degree $n$ with exactly one complex place. The conjectural higher elliptic units in \cite{thirdpaper} are given by evaluations of the $(n-2)$-cocycles $\psi_{n,a}^{(N)}$ against some $(n-2)$-cycles on these groups of totally positive units.

\bibliographystyle{alpha}
\bibliography{bibliographie}{}

\end{document}